\documentclass[10pt,a4paper]{article}
\pdfoutput=1
\usepackage{amsmath,amsthm,amssymb}
\usepackage
   {graphicx}
\usepackage{url} 
\title{Approximating optimization problems over convex functions}
\author{N\'{e}stor E. Aguilera \and Pedro Morin%
   \thanks{This work is supported in part by CONICET (Argentina) Grant PIP 5810.}
   }

\theoremstyle{plain}
   \newtheorem{theorem}{Theorem}[section]
   \newtheorem{corollary}[theorem]{Corollary}
   \newtheorem{lemma}[theorem]{Lemma}
\theoremstyle{definition}
   
\theoremstyle{remark}
   \newtheorem{problem}[theorem]{Problem}
   \newtheorem{remark}[theorem]{Remark}
   \newtheorem{example}[theorem]{Example}
\newcommand{\ENDsymbol}{$\Diamond$}
\newcommand{\END}{~\hfill\ENDsymbol}
\numberwithin{equation}{section}
\numberwithin{figure}{section}
\numberwithin{table}{section}
\newcommand{\NN}{\mathbb N}
\newcommand{\ZZ}{\mathbb Z}
\newcommand{\RR}{\mathbb R}
\newcommand{\Conv}{\mathcal{C}}

\newcommand{\mesh}{\mathcal{M}}
\newcommand{\Q}{Q}
\newcommand{\Seq}{\mathcal{S}}
\newcommand{\U}{\mathcal{U}}
\newcommand{\dd}{\mathrm{d}}
\DeclareMathOperator{\dist}{dist}
\DeclareMathOperator{\diff}{\Delta}
\newcommand{\grad}{\nabla}
\newcommand{\pderiv}[2][{}]{\partial^{#1}_{#2}}
\newcommand{\abs}[1]{{\lvert#1\rvert}}
\newcommand{\bdy}{\partial}
\newcommand{\norm}[1]{{\lVert#1\rVert}}
\newcommand{\normast}[1]{{\lVert#1\rVert}^\ast}
\newcommand{\vv}[1]{\mathbf{#1}} 
\newcommand{\vvv}[1]{#1} 
\newcommand{\e}{{\vv{e}}}
\newcommand{\interior}[1]{\mathring{#1}}
\newcommand{\Ho}{\overline{H}}
\newcommand{\Hw}{\widetilde{H}}
\newcommand{\aut}[1]{\textsc{#1}}
\newcommand{\tita}[1]{#1}           
\newcommand{\titb}[1]{\emph{#1}}    
\newcommand{\titj}[1]{\emph{#1}}    
\begin{document}
\maketitle
\thispagestyle{empty}

\begin{abstract}

Many problems of theoretical and practical interest involve finding an optimum over a family of convex functions.
For instance, finding the projection on the convex functions in $H^k(\Omega)$, and optimizing functionals arising from some problems in economics.

In the continuous setting and assuming smoothness, the convexity constraints may be given locally by asking the Hessian matrix to be positive semidefinite, but in making discrete approximations two difficulties arise: the continuous solutions may be not smooth, and functions with positive semidefinite discrete Hessian need not be convex in a discrete sense.

Previous work has concentrated on non-local descriptions of convexity, making the number of constraints to grow super-linearly with the number of nodes even in dimension 2, and these descriptions are very difficult to extend to higher dimensions.

In this paper we propose a finite difference approximation using positive semidefinite programs and discrete Hessians, and prove convergence under very general conditions, even when the continuous solution is not smooth, working on any dimension, and requiring a linear number of constraints in the number of nodes.

Using semidefinite programming codes, we show concrete examples of approximations to problems in two and three dimensions.

\end{abstract}

\section{Introduction}
\label{sec:intro}

Convex and concave functions appear naturally in many branches of science such as biology (growth), medicine (dose-response), or economics (utility, production or costs), spurring in turn the interest of other areas, for example, statistics.
It is no surprise, then, that many problems of theoretical and practical interest involve the optimization of a functional over a family of convex functions.

A particularly important case is when the functional to be minimized is the norm of some normed function space $V$, i.e., given $f\in V$ find
\[ \min_{u\in \Conv}\, \norm{u - f}_V, \]
where $\Conv$ is a family of convex functions in $V$.
Typical choices are the spaces $L^2(\Omega)$ or $H^1(\Omega)$, and, in general, the $W^{k,p}(\Omega)$ spaces, where $\Omega$ is a convex domain in $\RR^d$.

Sometimes the convexity is a reasonable shape assumption on the model, which could be replaced by or added to other shape constraints such as radial symmetry, harmonicity or upper and lower bounds.
This is the case, for example, of Newton's problem of minimal resistance (see, e.g., \cite{BG97,BFK95,LP99,LP01}).

More surprisingly perhaps, the convexity constraint may be a consequence of the model, as in the design of some mechanisms in economics~\cite{MV07,RC98}.
Actually, our interest in the subject arose from one of these problems, which we will call the \emph{monopolist problem}, and is described in some detail in section~\ref{sec:economics}.
In this problem we wish to find
\begin{equation}\label{equ:econo}
   \max_{u\in\Conv}
      \int_\Q \bigl( \grad u(x) \cdot x - u(x) \bigr) \,f(x)\, \dd{x},
\end{equation}
where
$\Q = [0,1]^d$,
$f$ is a non-negative probability density function over $\Q$,
and $\Conv$ is a family of convex functions on $\Q$ with some further properties.

The monopolist problem is numerically very challenging since, unlike the problems coming from physics, the space dimension $d$ may be much higher than $2$ or $3$.

There is a big qualitative jump going from one to more dimensions when dealing with convexity constraints.
This can be appreciated readily by looking at the statistics literature, where the one dimensional convex regression or density problem has been considered both theoretically and numerically, and using different measuring functionals (see, e.g., \cite{ BD07,GJW01,HP76, Me06}), but very little has been done numerically in two or more dimensions.
In this regard, it is worth mentioning the work by Shih, Chen and Kim~\cite{SCK06}, who use appropriate splines in the MARS (Multivariate Adaptive Regression Splines) algorithm.

One of the main difficulties in obtaining discrete approximations to optimization problems over convex functions lies in giving a local and finite description of convex functions if $d > 1$.
Though this could be done for smooth functions of continuous variables by asking the Hessian matrix to be positive semidefinite at all points, there is no similar characterization for discrete functions on meshes.
As a matter of fact, we show in examples~\ref{example:hessnoconv} and~\ref{example:convnohess} that this cannot be done easily.

It goes without saying that there is a lot of work done on convex functions of continuous variables, as exemplified by the classical book by Rockafellar~\cite{rockafellar}.
Traditionally, this work leans more on properties satisfied by convex functions rather than on properties that imply convexity.
In particular, very little has been done on local properties guaranteeing convexity.

For discrete variables, there is a rather large work done by the discrete mathematics community for fixed lattices, and a number of definitions for discrete convex functions have been proposed (see, e.g., \cite{murota}).
Again, usually these definitions are non-local.

As far as we know, there exists very few literature on optimization problems on convex functions in dimensions two or more, either theoretically or numerically.
Besides those already mentioned (and the references therein), let us cite three more which are prominent in the context of our work.

\begin{itemize}
\item
Carlier and Lachand-Robert~\cite{carlier-cpam} obtained the $C^1$ regularity of a variant of the monopolist problem, substituting the functional in~\eqref{equ:econo} for
\[
   \int_{\Omega} \Bigl(- \frac{1}{2}\,\abs{\grad u}^2
      + x\cdot \grad u - u(x) \Bigr)\,f(x)\,\dd{x},
\]
under some restrictions on the domain $\Omega$ and the density $f$.
They obtained also $C^1$ regularity for convex minimizers of functionals of the form
\[
   \int_{\Omega} \bigl(
      A(x,u)\grad u(x)\cdot \grad u (x) + f(x) u(x)\bigr)\,\dd{x}.
\]

\item
Carlier, Lachand-Robert and Maury~\cite{carlier} proposed a numerical scheme for minimizers of functionals of the type
\[ \int_\Omega j(x, u(x), \grad u(x))\,\dd{x}, \]
on closed convex subsets of convex functions of $H^1(\Omega)$ or $L^2(\Omega)$, where $j$ is a quadratic function of $u$ and $\grad u$, and $\Omega\subset\RR^2$.
We will discuss their approach in the next section (after the inequalities~\eqref{equ:discrete:conv}).

As Carlier et al.~\cite{carlier} point out, their work encompasses
the problem of finding
\[
   \min\int_{\Omega} \abs{u - f}^2\,\dd{x}
   \qquad\text{subject to}\quad
   u\in L^2(\Omega),\ \text{$u$ convex},\ u\le f,
\]
for given $f\in L^2(\Omega)$, i.e., a $L^2$-norm projection, and this problem is equivalent to that of finding the convex envelope $f^{**}$ of $f$.
Thus, minimizing over convex functions and finding the convex envelope of a function are two quite related tasks.

\item
Being one of the main problems in computational geometry, there are a number of well established codes for finding the convex hull of a set of points in $\RR^d$, which are very efficient in low dimensions.
Hence, it is natural to try to use these codes to approximate optimization problems on convex functions, an approach which Lachand-Robert and Oudet~\cite{LO05} applied to several problems.

It would be very interesting to see whether these ideas could be carried over, since the convex hull codes are quite fast in low dimensions.

\end{itemize}

Our approach, based on semidefinite programming, takes a different direction from those mentioned.

Let us recall that a semidefinite program is an optimization problem of the form
\begin{equation}\label{equ:psd}
\begin{gathered}
   \min c\cdot x\\
   \begin{array}{rc}
   \text{subject to}&\\
   & x_1 A_1 + \dots + x_n A_n - A_0 \succeq 0,\\
   &x\in\RR^n,
   \end{array}
\end{gathered}
\end{equation}
where $c\in\RR^n$, $A_0,A_1,\dots,A_n$ are symmetric $m\times m$ matrices, and $A\succeq 0$ indicates that the symmetric matrix $A$ is positive semidefinite.
By letting the matrices $A_i$ be diagonal, we see that the program~\eqref{equ:psd} is a generalization of linear programming (and includes it strictly).
Thus, in a semidefinite program the constraints can be a mixture of linear inequalities and positive semidefinite requirements.

In this paper we give a theoretical framework for approximating many optimization problems on convex functions using a finite differences scheme which imposes a positive semidefinite constraint on a discretization of the Hessian matrix.

Although not linear, our approach seems very natural and has many advantages.
Being of a local nature, the number of constraints grows only linearly with the number of nodes, and it works for any dimension of the underlying space.
Notwithstanding the already mentioned counterexamples (\ref{example:hessnoconv} and~\ref{example:convnohess}) of the relation between convexity and positive semidefinite Hessian, we will show that for many problems we obtain convergence to a continuous optimum.
This convergence holds even for non-smooth optima, as those arising when projecting in the $L^\infty$ norm or in some problems of the type given by equation~\eqref{equ:econo}.

In practice, our definition can be used to advantage by using the existing efficient semidefinite codes and, furthermore, it is very simple to program in higher dimensions.

As a final remark, we think that our approach might be a first step to deal with other problems where convexity is a consequence of the model, such as transportation problems where the cost is quadratic, i.e., for solving the Monge-Amp\`ere equations, a prime example of fully non-linear differential equations.

This article is organized as follows.

In section~\ref{sec:review} we summarize some techniques that could be used to deal with the numerical approximation of optimal convex functions, showing their strengths and weaknesses, and introduce the discrete Hessian.

In section~\ref{sec:convex:cont} we give different characterizations of smooth convex functions of continuous variables, which allow us to extend the definition of the Hessian to $C^1$ functions in terms of averages.

The heart of the paper is section~\ref{sec:approxs}: in order to have a good definition of discrete convexity, we need to show that any convex function of continuous variables may be approximated by its discrete counterparts, and conversely, that a converging sequence (in a suitable norm) of discrete convex functions will do so to a convex function of continuous variables.

We show in this section that the approximation of the discrete functions to the limit $u$ is uniform in $u$ and its derivatives if $u$ is smooth, and uniform in $u$ over compact subsets of $\Omega$, if $u$ is merely convex and bounded (and hence continuous).
Although our definition is stated using the Hessian, if a sequence of discrete functions converges to a function which is not $C^2$ or even $C^1$, still convexity of the limit function is guaranteed.

In section~\ref{sec:functionals} we pose a general structure for optimization problems on convex functions which fits the results of the previous sections, and may be applied to several important problems.
In addition, error bounds may be obtained if the continuous optimal solution is smooth.

The results of sections~\ref{sec:convex:cont}--\ref{sec:functionals} are somewhat independent of semidefinite programming.
However, as explained in section~\ref{sec:review}, for many problems of interest the functional and constraints involved may be expressed as a semidefinite program, and we present in section~\ref{sec:numerical} numerical examples showing that the current codes make it feasible to use our scheme on them.

We begin section~\ref{sec:numerical} by considering the monopolist problem in two and three dimensions, comparing our results to the analytical solution when $f = 1$ in~\eqref{equ:econo}.
We continue by showing some norm projections, exploring the behavior under different functionals of an example given by Carlier, Lachand-Robert and Maury~\cite{carlier}.
Moreover, we exhibit an explicit example in which the $L^\infty$ projection is not unique.
We finish the section on numerical examples by showing how our scheme could be used to fit a (discrete) convex function to noisy data given on the nodes of a regular mesh, using the $L^2$ norm (i.e., a least squares approach), and also the $L^1$ and $L^\infty$ norms, which are not often seen.

\section{Dealing with convex functions}
\label{sec:review}

In this section we present some ideas and techniques that could be used to approximate numerically an optimal convex function.
For the sake of simplicity, we will work on the unit $d$-dimensional cube in $\RR^d$, $\Q = [0,1]^d$.

If $d = 1$ and
\begin{equation}\label{equ:1d:mesh}
   x_0 = 0 < x_1 < \dots < x_n = 1,
\end{equation}
is a discretization of $[0,1]$, then a convex function $u$ defined on $[0,1]$ satisfies the inequalities
\begin{equation}\label{equ:1d:convex}
   \frac{u(x_{i}) - u(x_{i-1})}{x_{i} - x_{i-1}}
   \le \frac{u(x_{i+1}) - u(x_i)}{x_{i+1} - x_i}
   \quad\text{for $i = 1,\dots,n-1$}.
\end{equation}

Conversely, if a function $u$ defined only in the mesh points $\{x_0,\dots,x_n\}$ satisfies the inequalities~\eqref{equ:1d:convex}, we can always extend it linearly in the intervals $[x_{i-1},x_i]$, $i = 1,\dots,n$, so that the resulting piecewise linear function will be convex on $[0,1]$.

Thus, the inequalities~\eqref{equ:1d:convex} may be used to define discrete convexity in one dimension.
It is clear that---via the piecewise linear extension---we can approximate any convex function in $[0,1]$.
On the other hand, as we show in lemma~\ref{lemma:1d:convergence:1}, under some adequate assumptions if a sequence of discrete convex functions (satisfying~\eqref{equ:1d:convex} for each corresponding subdivision) converges pointwise, then the limit defines a unique convex function in $[0,1]$.

Let $h$ be given, $0 < h < 1/2$, and suppose $\mesh = \{x_0,\dots,x_n\}$ is a subdivision of $[0,1]$ satisfying~\eqref{equ:1d:mesh} and such that
\begin{equation}\label{equ:1d:size}
   \max_{1\le i\le n} \abs{x_i - x_{i-1}} \le h.
\end{equation}
Let $u$ be defined on $\mesh$ and satisfy~\eqref{equ:1d:convex}. For a given $\varepsilon$, $h < \varepsilon < 1/2$, let us consider
\[
   \underline{x} = \max\, \{x\in\mesh : x \le \varepsilon\}
   \quad\text{and}\quad
   \overline{x} = \min\, \{x\in\mesh : x \ge 1 - \varepsilon\}.
\]
From~\eqref{equ:1d:convex} we see that
\[
   -\frac{|u(\underline{x})-u(0)|}{\varepsilon}
   \le \frac{u(\underline{x})-u(0)}{\underline{x}}
   \le \frac{u(x_{i+1}) - u(x_i)}{x_{i+1} - x_i}
   \le \frac{u(1) - u(\overline{x})}{1 - \overline{x}}
   \le \frac{|u(1) - u(\overline{x})|}{\varepsilon}
\]
for all $x_i,x_{i+1}\in\mesh$ such that $\underline{x} \le x_{i}$, $x_{i+1} \le \overline{x}$.
Therefore, if
\[ \abs{u(x_i)} \le M \quad\text{for all $x_i\in\mesh$,} \]
we have
\begin{equation}\label{equ:conv:bound:deriv}
   \underline{C}
   \le \frac{u(x_{i+1}) - u(x_i)}{x_{i+1} - x_i}
   \le \overline{C}
   \quad
   \text{for all $\varepsilon \le x_{i}$, $x_{i+1} \le 1 - \varepsilon$,}
\end{equation}
where $\underline{C}$ and $\overline{C}$ are constants depending only on $\varepsilon$ and $M$ (but not on $u$ or $\mesh$).
In other words, if $u$ is extended as a piecewise linear function to all of $[0,1]$, the resulting function is uniformly Lipschitz on compact subsets of $(0,1)$.

Hence,
{\lemma\label{lemma:1d:convergence:1}
Let $(\mesh_h)_h$ be a sequence of meshes in $[0,1]$, $h\downarrow 0$, satisfying~\eqref{equ:1d:size} for each $h$, and such that $\mesh_h \subset \mesh_{h'}$ for $h > h'$.
Suppose $u_h$ is defined in $\mesh_h$ for each $h$, and the sequence $(u_h)_h$ is such that for all $x\in \cup_h \mesh_h$, $\lim_h u_h(x)$ exists and is finite, and let
\[ u(x) = \lim_h u_h(x) \quad\text{for $x\in\cup_h \mesh_h$.} \]

Then $u$ may be extended to all of $[0,1]$ as a continuous convex function in a unique way.
Moreover, the convergence of $u_h$ (extended piecewise linearly) to $u$ is uniform in compact subsets of $(0,1)$.
\endlemma}

Also, by the Arzela-Ascoli theorem,

{\lemma\label{lemma:1d:convergence:2}
Let $(\mesh_h)_h$ be a sequence of meshes as in the previous lemma.
Suppose $M > 0$ is given, and that for each $h$ a function $u_h$ is defined in $\mesh_h$ so that
\[
   \abs{u_h(x)} \le M
   \quad\text{for all $x\in \mesh_h$ and all $h$,}
\]
and assume $u_h$ is extended to all of $[0,1]$ piecewise linearly.

Then, there exists a subsequence of $(u_{h'})$ of $(u_h)$ converging to a continuous convex function $u$ defined on $[0,1]$.
Moreover, the convergence of the subsequence $(u_{h'})$ to $u$ is uniform in compact subsets of $(0,1)$.
\endlemma}

Thus, from a theoretical point of view, the discrete functions satisfying~\eqref{equ:1d:convex} are well understood.

From the numerical point of view, when used in a discretization of an optimization problem, the constraints coming from the inequalities~\eqref{equ:1d:convex} are linear, and solving the resulting discrete optimization problem with them is usually not much harder than solving it without them.

Hence, the case $d = 1$ poses no major trouble.

For $d > 1$ it will be more convenient to work only with regular meshes (or grids) on $\Q = [0,1]^d$.
Thus, for fixed $h > 0$ ($h = 1/n$ for some $n\in\NN$),
the mesh $\mesh_h$ will consist of all points $x\in\RR^d\cap \Q$ of the form $x = h z$ with $z\in\ZZ^d$.
Denoting by \(\interior{\Q} = \Q\setminus \bdy\Q\) the interior of $\Q$, we set
\[
   \interior{\mesh}_h = \mesh_h\cap \interior{\Q},
   \qquad
   \bdy\mesh_h = \mesh_h\cap \bdy\Q,
\]
and denote by $\U_h$ the set of real valued functions defined on $\mesh_h$.

A first simple idea to extend the inequalities~\eqref{equ:1d:convex} to more dimensions, is to consider the set of functions $u_h\in\U_h$ satisfying the convexity constraints
\[
   u_h(x + h\e_i) + u_h(x - h\e_i) \ge 2\,u_h(x)
   \quad\text{for all $x\in\interior{\mesh}_h$ and $i = 1,\dots,d$},
\]
where $\e_i$ denotes the $i$-th vector of the canonical basis of $\RR^d$.

This set of discrete functions is very appealing because the convexity is modelled by using only $O$(number of mesh points) linear constraints, as in the case $d = 1$.
As we have seen, this approach is exact in one dimension, and gives satisfactory results in many cases in more dimensions (in particular for some of the specific examples treated later in section~\ref{sec:economics}), but there is no guarantee of convexity in the limit function if $d > 1$.
In fact, by taking the interpolant, with this set we can certainly approximate any function of continuous variables which is convex, but we can also approximate other functions.
For example, we may approximate $u(x_1,x_2) = x_1 x_2$ which is linear---and thus convex---in each coordinate direction, but not convex in $\Q\subset\RR^2$.

This shows that the definition of discrete convexity should be done with more care: when $d > 1$, we have to take into account every possible direction.

It is reasonable, then, to say that a discrete function is convex if
\begin{equation}\label{equ:discrete:conv}
   u_h(x + y) + u_h(x - y) \ge 2u_h(x)
   \quad\text{for all $x,y$ such that $x\pm y\in\mesh_h$.}
\end{equation}
Since as $h$ goes to $0$ the possible directions become dense, it is possible (under some conditions) to regain convexity in the limit as we had for one dimension.

In a way, this is the approach followed by Carlier, Lachand-Robert and Maury~\cite{carlier}.
In a two dimensional setting, they consider discrete convex functions which are restrictions to a mesh of convex functions of continuous variables, and show that this definition is equivalent to an intrinsic one, stated only in terms of the value of the function at the grid points, similar to~\eqref{equ:discrete:conv}.
The problem with this description is that it is non-local, and the number of constraints needed in two dimensions (after pruning) reportedly grows approximately as $N^{1.8}$, where $N$ is the number of nodes in the grid.
Moreover, this approach is very difficult to extend to higher dimensions.

In order to keep the definition of discrete convexity local, another possibility is to consider discretizations of the Hessian matrix, as we do in this work.

For $u\in\U_h$ and $x\in\mesh_h$, we define the (forward) first order finite differences by
\[
   \diff_{h,i} u (x) = \frac{u(x + h\,\e_i) - u(x)}{h},
   \qquad i = 1,\dots,d,
\]
and the second order finite differences by
\[
   \diff_{h,ii}^2 u(x)
      = \frac{u(x + h\,\e_i)
         - 2\,u(x) + u(x - h\,\e_i)}{h^2},
   \qquad i = 1,\dots,d,
\]
or, if $i\ne j$, by
\begin{multline*}
   \diff_{h,ij}^2 u (x) = \frac{1}{4\,h^2}\,
   \bigl(
        u(x + h\,\e_i + h\,\e_j) - u(x - h\,\e_i + h\,\e_j) \\
      - u(x + h\,\e_i - h\,\e_j) + u(x - h\,\e_i - h\,\e_j)
   \bigr).
\end{multline*}

Clearly, not all of these finite differences are defined for points in $\bdy\mesh_h$.
Therefore, when mentioning $\diff_{h,i} u(x)$ or $\diff^2_{h,ij} u(x)$ for $x\in\mesh_h$, we will implicitly assume that $x$ and $i$ (and eventually $j$) are such that the corresponding finite difference is well defined at $x$.

Finally, we define the discrete Hessian of $u\in\U_h$ at $x\in\interior{\mesh}_h$, as the symmetric matrix $H_h u(x)\in\RR^{d\times d}$ whose $i,j$ entry is $\diff_{h,ij}^2 u(x)$.

We are faced now with two issues:

{\enumerate
\item
How can a discrete optimization problem with positive semidefinite discrete Hessian constraints be solved?

\item
Once we found a discrete solution, how does it approximate the continuous solution?
\endenumerate}

As mentioned in the introduction, for the first issue we will use the semidefinite programming model~\eqref{equ:psd}, where the objective is linear and the constraints are positive semidefinite.
Thus, if the objective of the continuous problem is not linear, we must rewrite it as a constraint.
However, it is worth mentioning that not every functional may be readily modelled as a positive semidefinite constraint, an example of such functionals being
\[ \int_\Omega \frac{1}{1 + \abs{\grad u}^2}\,\dd{x}, \]
arising in Newton's problem of minimal resistance.

Although this method may not be used for \emph{all} functionals, it can be used in \emph{many} cases as we now illustrate.
The interested reader is referred to the article by Vandenberghe and Boyd~\cite{vandenberghe} for other possibilities.

In what follows we associate with each mesh node $P_k$, $k = 1,\dots,N$,
the unknown value $u_k$,
the discrete Hessian $H_h(P_k)$,
and a $d$-dimensional cube $Q_k$ centered at $P_k$ of measure $\abs{Q_k}$, so that $\sum_k \abs{Q_k} = \abs{\Q}$.

{\description
\item[$L^1$ projection.]
The program
\[
   \min\, \int_{Q} \abs{u - f} \,\dd{x} =
   \min\, \sum_k \int_{Q_k} \abs{u - f} \,\dd{x},
\]
subject to $u$ convex, is modeled by adding $N$ more unknowns $t_1,\dots,t_N$ (for a total of $2N$) as
\[
\begin{gathered}
   \min\, \sum_k t_k\\
   \begin{array}{rc}
   \text{subject to}\\
      &\begin{aligned}
       \int_{Q_k} \abs{u - f} \,\dd{x} &\le t_k,
         && k = 1,\dots,N,\\
       H_h(P_k) &\succeq 0, && k = 1,\dots,N.
      \end{aligned}
   \end{array}
\end{gathered}
\]

In turn, the constraints of the form
\[ \int_{Q_k} \abs{u - f}\,\dd{x} \le t_k, \]
are lumped as $\abs{u_k - f(P_k)}\,\abs{Q_k} \le t_k$, and finally modeled as the linear inequalities
\[\begin{split}
   -u_k \abs{Q_k} + t_k &\ge -f(P_k)\,\abs{Q_k},\\
    u_k \abs{Q_k} + t_k &\ge  f(P_k)\,\abs{Q_k}.
\end{split}\]

\item[$L^2$ projection.]
This is analogous to the $L^1$ projection:
\[
\begin{gathered}
   \min\, \sum_k t_k\\
   \begin{array}{rc}
   \text{subject to}\\
   & \begin{aligned}
      (u_k - f_k)^2\,\abs{Q_k} &\le t_k, && k = 1,\dots,N\\
      H_h(P_k) &\succeq 0, && k = 1,\dots,N,
      \end{aligned}
   \end{array}
\end{gathered}
\]
and the constraints of the form \((u_k - f_k)^2\,\abs{Q_k} \le t_k\) are written in the form
\[
   \begin{bmatrix}
      t_k & (x_k - f_k)\,\sqrt{\abs{Q_k}}\, \\
      (x_k - f_k)\,\sqrt{\abs{Q_k}} & 1
   \end{bmatrix} \succeq 0.
\]

The $H^1$ projection is handled in a similar way.

\item[$L^\infty$ projection.]
This is analogous to the other projections, except that we only need to add one more variable $t$ (instead of $N$), obtaining the discrete program
\[\begin{gathered}
   \min\, t\\
   \begin{array}{rc}
   \text{subject to}\\
   & \begin{aligned}
      -u_k + t &\ge -f(P_k), && k = 1,\dots, N,\\
      u_k + t &\ge  f(P_k), && k = 1,\dots,N,\\
      H_h(P_k) &\succeq 0, && k = 1,\dots,N.
      \end{aligned}
   \end{array}
\end{gathered}\]
\enddescription}

The next issue is to see how well the discrete solutions of the semidefinite program approximate the continuous convex solution.

Our first hope is that $H_h u(x) \succeq 0$ will be equivalent to some version of convexity of discrete functions, for example to the one given by the inequalities~\eqref{equ:discrete:conv}.
The next examples show that this is not true.
The first one shows that a non-convex function may have a positive semidefinite discrete Hessian, and the second one shows that the discrete Hessian may not be positive semidefinite for some convex functions.

{\example\label{example:hessnoconv}
Let us consider $d = 2$, and $u$ defined on the 2-dimensional $2\times 2$ grid $\mesh_{1/2}$ by
\[
\begin{aligned}
   u(0,1) &= 1, & u(1/2,1) &= 1, & u(1,1) &= 1, \\
   u(0,1/2) &= 1/2, & u(1/2,1/2) &= 5/8, & u(1,1/2) &= 1, \\
   u(0,0) &= 0, & u(1/2,0) &= 1/2, & u(1,0) &= 1.
\end{aligned}
\]

We may check that (with $h = 1/2$),
\[
   H_h u(1/2,1/2) =
   \begin{bmatrix}
   1 & -1\\
   -1 & 1
   \end{bmatrix},
\]
whose eigenvalues are $2$ and $0$ (with eigenvectors \((1,-1)\) and \((1,1)\)), and is thus positive semidefinite.

However, the restriction to the diagonal $x_1 = x_2$ is not convex:
\[
   \frac{u(0,0) + u(1,1)}{2} = \frac{1}{2} < \frac{5}{8} = u(1/2,1/2).
   \tag*{\ENDsymbol}
\]
\endexample}

\smallskip
{\example\label{example:convnohess}
Consider the same grid of the previous example and $u$ defined by
\[
\begin{aligned}
      u(0,1) &= 8/15,
      & u(1/2,1) &= 8/15,
      & u(1,1) &= 1,
      \\
      u(0,1/2) &= 1/30,
      & u(1/2,1/2) &= 1/2,
      & u(1,1/2) &= 1,
      \\
      u(0,0) &= 0,
      & u(1/2,0) &= 1/2,
      & u(1,0) &= 1.
\end{aligned}
\]

We see that $u$ is the restriction to $\mesh_h$ of the convex function
\[
   \max\,
      \Bigl\{x_1,
      \frac{14 x_1 + x_2}{15},
      x_2 - \frac{7}{15}
      \Bigr\},
\]
but
\[
   H_h u(1/2,1/2) = \frac{1}{15}\,
      \begin{bmatrix} 2 & -23\\ -23 & 2 \end{bmatrix},
\]
which has eigenvalues $5/3$ and $-7/5$, and hence is not positive semidefinite.\END
\endexample}

\smallskip

Notwithstanding these examples, in section~\ref{sec:approxs} we will show that, under suitable conditions, by asking for positive semidefinite discrete Hessians we may conveniently approximate continuous convex functions, and obtain variants of lemmas~\ref{lemma:1d:convergence:1} and~\ref{lemma:1d:convergence:2}.

However, we will need some previous results which are tackled in the following section.

\section{Convex functions of continuous variables}
\label{sec:convex:cont}

In this section we deal with variants of the Hessian matrix, initially defined for $C^2$ functions, which allow us to characterize convexity for $C^1$ functions locally.

Let us start by introducing some more precise notation:

{\itemize
\item
The distance from a point $x$ to a set $A$ will be denoted by $\dist(x,A)$.

\item
If $x\in\RR^d$ and $h > 0$, \(\Q_h(x) = \{y\in\RR^d : \abs{y_i - x_i} \le h \text{\ for all $i = 1,\dots,d$}\}.\)

\item
The gradient of $u$ is denoted by
\(\grad u = (\pderiv{1} u, \dots, \pderiv{d} u)\),
and we write $\pderiv[2]{ij}$ for $\pderiv{i} \pderiv{j}$.

\item
If $\Omega\subset\RR^d$ is an open set, we denote by $C^k(\Omega)$ the set of functions having all derivatives up to order $k$ continuous on $\Omega$, and by $C^k(\overline{\Omega})$ the set of functions having all derivatives up to order $k$ uniformly continuous in $\Omega$.
If $\Omega$ is omitted, $C^k = C^k(\Q)$.
\enditemize}

Let $\varphi$ be a non-negative, real valued function in $C_0^\infty$, vanishing outside \(\{x\in\RR^d : \abs{x} < 1\}\), and such that
\( \int_{\RR^d} \varphi(x) \,\dd{x} = 1, \)
and for $\varepsilon > 0$ define
\(
   \varphi_\varepsilon (x)
   = \varepsilon^{-d}\, \varphi(\varepsilon^{-1} x).
\)
The function
\begin{equation}\label{equ:approx:epsilon}
   u^\varepsilon (x)
   = u\ast \varphi_\varepsilon (x)
   = \int_{\RR^d} u(x - y)\,\varphi_\varepsilon(y) \,\dd{y},
\end{equation}
defined for functions $u$ and points $x$ whenever the right hand side makes sense, has many interesting well known properties, and we state some of them without proof, referring the reader to, e.g., the book by Ziemer~\cite[Theorem~1.6.1 and Remark~1.6.2]{ziemer}.

{\theorem\label{theorem:ziemer}
\mbox{}\par
{\enumerate
\item
If $u\in L^1_\text{loc} (\RR^d)$, then for every $\varepsilon > 0$, $u^\varepsilon \in C^\infty(\RR^d)$ and
\( \pderiv[\alpha]{} u^\varepsilon = (\pderiv[\alpha]{} \varphi_\varepsilon)\ast u, \)
for every multi-index $\alpha$, where for $\alpha = (\alpha_1,\dots,\alpha_d)$, $\pderiv[\alpha]{} = \pderiv[\alpha_1]{1}\dots\pderiv[\alpha_d]{d}$.

\item
If $u\in L^1(\Omega)$ then $u^\varepsilon (x)$ is defined for $x\in\Omega$ and $\dist(x,\bdy\Omega) > \varepsilon$.

\item
\(\lim_{\varepsilon \to 0} u^\varepsilon (x) = u(x)\)
whenever $x$ is a Lebesgue point of $u$.

\item
If $u\in C(\Omega)$, and $F$ is a compact subset of $\Omega$, then $u^\varepsilon$ converges uniformly to $u$ on $F$ as $\varepsilon\to 0$.
\endenumerate}
\endtheorem}

The functions $u^\varepsilon$, also called regularizations or mollifiers of $u$, will make us work often with the sets
\[ A_\varepsilon = \{x \in\Q : \dist(x,\bdy\Q) > \varepsilon\}, \]
where $\varepsilon\in\RR$, $0 < \varepsilon < 1$.

Our first result is a simple characterization of continuous convex functions using the regularizations $u^\varepsilon$, and we omit its proof.

{\lemma\label{lemma:cont:1}
Suppose $u\in C$ and $u^\varepsilon$ is defined as in~\eqref{equ:approx:epsilon}.
We have,
{\enumerate
\item
If $u$ is convex, then $u^\varepsilon$ is convex in $A_\varepsilon$ for all $\varepsilon > 0$.

\item
Conversely, if for a sequence of $\varepsilon$'s converging to $0$, $u^\varepsilon$ is convex in $A_\varepsilon$, then $u$ is convex.
\endenumerate}
\endlemma}

The Hessian of $u\in C^2$ at $x$ is the symmetric matrix $H u(x)$ whose $ij$ entry is $\pderiv[2]{ij} u(x)$, and convex functions in $C^2$ are characterized by the positive semidefinite condition
\[ Hu(x)\succeq 0 \quad\text{for all $x\in\interior{\Q}$}. \]

Also, for $u\in C^2$, and $\delta > 0$, the average
\begin{equation}\label{equ:cont:def}
   \Ho_{\delta} u(x)
   = \frac{1}{(2 \delta)^d}\,\int_{\Q_{\delta}(x)} \,Hu(y)\,\dd{y},
\end{equation}
defined for $x\in A_\delta$, converges to $Hu(x)$ for all $x\in\interior{\Q}$, as $\delta\downarrow 0$.
Thus the condition
\begin{equation}\label{equ:cont:cond1}
   \Ho_{\delta} u(x) \succeq 0
   \quad
   \text{for all $x$ and $\delta$ such that $x\in A_\delta$},
\end{equation}
implies the convexity of $u$.
Actually, since the average of positive semidefinite matrices is positive semidefinite, \eqref{equ:cont:cond1} is equivalent to convexity for $u\in C^{2}$.

The definitions of $Hu$ and $\Ho_{\delta} u$ involve second order derivatives which are not defined if $u$ is not smooth enough.
However, we may express $\Ho_{\delta} u$ merely in terms of $u$ and $\grad u$ integrating by parts, and the resulting formula will make sense for $u\in C^{1}$.
Thus, for $u\in C^1$ and $x\in A_\delta$ let us define $\Hw_\delta u(x)$ as the symmetric ${d\times d}$ matrix whose diagonal entries are
\begin{subequations}\label{equs:cont:2}
\begin{equation}\label{equ:cont:2:ii}
   \bigl(\Hw_\delta u(x)\bigr)_{\!ii}
   = (2\delta)^{-d}\, \Bigl(
       \int_{F^{+}_{\delta,i}(x)} \pderiv{i} u(y)\,\dd{y}
     - \int_{F^{-}_{\delta,i}(x)} \pderiv{i} u(y)\,\dd{y}
     \Bigr),
\end{equation}
and, if $i\ne j$,
\begin{multline}
\label{equ:cont:2:ij}
   \bigl(\Hw_\delta u(x)\bigr)_{ij}
   = (2\delta)^{-d}\,
      \Bigl(
   \int_{F^{+}_{\delta,i}(x)\cap F^{+}_{\delta,j}(x)} u(y)\,\dd{y}
      - \int_{F^{+}_{\delta,i}(x)\cap F^{-}_{\delta,j}(x)} u(y)\,\dd{y}
   \\[4pt]
    -
    \int_{F^{-}_{\delta,i}(x)\cap F^{+}_{\delta,j}(x)} u(y)\,\dd{y}
      + \int_{F^{-}_{\delta,i}(x)\cap F^{-}_{\delta,j}(x)} u(y)\,\dd{y}
      \Bigr),
\end{multline}
where $F^{\pm}_{\delta,i}(x)$ denotes the $d-1$ dimensional face of the cube $\Q_\delta(x)$ having outward normal $\pm \e_i$.
\end{subequations}

Of course, since the equations~\eqref{equs:cont:2} were obtained integrating by parts, we have:

{\lemma\label{lemma:cont:2}
If $u\in C^2$ then $\Ho_{\delta} u = \Hw_\delta u$.
\endlemma}

As we show now, the extension $\Hw_\delta$ of $\Ho_{\delta}$ to functions in $C^1$, still gives a local characterization of convexity.

{\theorem\label{theorem:cont:1}
The following are valid if $u\in C^1$:

{\enumerate
\item
If $u$ is convex, then
\( \Hw_\delta u (x) \succeq 0\)
for all $\delta > 0$ and $x\in A_\delta$.

\item
If for a sequence $\delta_n\downarrow 0$,
\( \Hw_{\!\delta_n} u (x) \succeq 0\)
for all $x\in A_{\delta_n}$, then $u$ is convex.
\endenumerate}
\endtheorem}

{\proof
Suppose $u\in C^1$, and consider the regularization $u^\varepsilon$ defined in~\eqref{equ:approx:epsilon}.
Since we can interchange derivatives and integrals with the convolution,
\begin{equation}\label{equ:pf1:1}
   \Hw_\delta u^\varepsilon (x)
   = \Hw_{\delta} (u \ast \varphi_\varepsilon) (x)
   = \bigl(\Hw_\delta u\bigr) \ast \varphi_\varepsilon (x)
   \quad\text{for $x\in A_{\delta + \varepsilon}$}.
\end{equation}

If $u$ is convex, by the first part of lemma~\ref{lemma:cont:1} we know that $u^\varepsilon$ is convex in $A_{\varepsilon}$, and since $u^\varepsilon\in C^\infty$, for $x\in A_{\delta + \varepsilon}$ we have
\( \Ho_{\delta} u^\varepsilon (x) \succeq 0 \).
Thus, by lemma~\ref{lemma:cont:2},
\[
   \Hw_\delta u^\varepsilon (x) \succeq 0
   \quad\text{for $x\in A_{\delta + \varepsilon}$}.
\]
Using theorem~\ref{theorem:ziemer}, since $u^\varepsilon$ and their derivatives converge uniformly to $u$ in compact subsets of $\interior{\Q}$, we must have $\Hw_\delta u(x)\succeq 0$ for $x$ in compact subsets of $A_\delta$, and therefore in the whole of $A_\delta$.

On the other hand, if $\Hw_\delta u (x) \succeq 0$ for $x\in A_\delta$, since an integral mean of positive semidefinite matrices is positive semidefinite, using~\eqref{equ:pf1:1} and lemma~\ref{lemma:cont:2}, we see that
\[
   \Ho_{\delta} u^\varepsilon (x) =
   \Hw_\delta u^\varepsilon (x)
   \succeq 0
   \quad\text{for $x\in A_{\varepsilon + \delta}$,}
\]
which implies that $u^\varepsilon$ is convex in $A_{\varepsilon + \delta}$.
Letting $\delta$ go to $0$ while keeping $\varepsilon$ fixed, we see that $u^\varepsilon$ is convex in $A_\varepsilon$, and by the second part of lemma~\ref{lemma:cont:1}, $u$ must be convex.
\endproof}

\section{Approximating convex functions, the discrete Hessian}
\label{sec:approxs}

There are two main issues when defining the set of discrete approximants to be used:

{\enumerate
\item we want it to be rich enough to approximate every convex function, and

\item we want this set to be not too large, to avoid convergence to non-convex functions.
\endenumerate}

The first point is very natural, and necessary to approximate the solution to the problem.
The second point might look artificial at first sight, but if not enforced, then we could be approximating a non-convex function.
This is the case, for instance, if we only require convexity along the coordinate axes, as we stated in section~\ref{sec:review}.

Our discrete approximations will be a subset of functions having positive semidefinite discrete Hessians, and the main purpose of this section is to address the two issues mentioned above.

We will show in theorem~\ref{theorem:discrete:approx} that, despite the examples~\ref{example:hessnoconv} and~\ref{example:convnohess}, the discrete Hessian $H_h$ may be used to obtain very good approximations to convex functions of continuous variables.

This is not surprising since we can approximate, say, $u\in C^3$, by discrete functions whose finite differences up to order $2$ converge to the derivatives up to order $2$ of $u$, and then add a small perturbation to bring up the eigenvalues of the discrete Hessians.
This is the main idea of the proof, and in its course, we will use the following part of the Hoffman-Wielandt theorem~\cite{hoffman-wielandt}.

{\theorem[Hoffman and Wielandt]
\label{theorem:hoffman}
There exists a positive constant $c_d$, depending only on the dimension $d$, such that if $A = [a_{ij}]$ and $B = [b_{ij}]$ are symmetric $d\times d$ matrices, and $\lambda$ and $\mu$ are their minimum eigenvalues, then
\[ \abs{\lambda - \mu} \le c_d\, \max_{ij}\,\abs{a_{ij} - b_{ij}}. \]
\endtheorem}

{\theorem\label{theorem:discrete:approx}
Let $u\in C^3$ be convex.
Then, for any $\varepsilon > 0$, there exists $h_0 = h_0(u,\varepsilon) > 0$ such that for any $h$, $0 < h < h_0$, there exists a function $u_h$ defined on $\mesh_h$ satisfying for all $x\in\mesh_h$,
\begin{multline}\label{equ:lemma:approx:3}
   \abs{u_h (x) - u(x)}
   + \sum_{i = 1}^d \abs{\diff_{h,i} u_h(x) - \pderiv{i} u(x)}\\
   + \sum_{1 \le i \le j \le d}
      \abs{\diff^2_{h,ij} u_h(x) - \pderiv[2]{ij} u(x)}
   < \varepsilon,
\end{multline}
and
\[ H_h u_h(x) \succeq 0 \quad\text{for all $x\in\interior{\mesh_h}$}. \]
\endtheorem}

Let us recall that in the inequality~\eqref{equ:lemma:approx:3}, for $x\in\bdy\mesh_h$ we consider only the finite differences that are defined at that $x$.

{\proof
If $u$ is convex and smooth, given $\varepsilon_1 > 0$ there exists $\delta_1 > 0$ such that for all $x\in\Q$ and $0 < h < \delta_1$,
\[
   \sum_{i = 1}^d \abs{\diff_{h,i} u(x) - \pderiv{i}u(x)}
   + \sum_{1 \le i \le j \le d}
      \abs{\diff^2_{h,ij} u(x) - \pderiv[2]{ij} u(x)}
   < \frac{\varepsilon_1}{c_d},
\]
where $c_d$ is the constant in theorem~\ref{theorem:hoffman}.
Hence, since $u$ is convex and therefore $Hu(x) \succeq 0$, the minimum eigenvalue of $H_h u(x)$ is uniformly bounded below by $-\varepsilon_1$.

Now consider the function
\[ g(x) = \frac{1}{2}\,\abs{x}^2, \]
for which, if $h$ small enough,
\[
   \abs{g(x)} \le d/2,
   \quad
   \abs{\diff_{h,i} g(x)} = \abs{x_i + h/2} \le 2,
   \quad\text{and}\quad
   H_h g(x) = H g(x) = I_d,
\]
where $I_d$ is the identity matrix in $\RR^{d\times d}$.

If $0 < h < \min\,\{\varepsilon_1, \delta_1\}$, the function
\[ u_h = u + \varepsilon_1 g, \]
defined on $\mesh_h$, satisfies
\begin{gather*}
   H_h u_h(x) \succeq 0 \quad\text{for $x\in\interior{\mesh_h}$},\\
   \abs{u_h(x) - u(x)} \le \frac{d\varepsilon_1}{2},
\end{gather*}
and
\begin{align*}
   &\hspace{-20pt}
   \sum_{i = 1}^d
      \abs{\diff_{h,i} u_h(x) - \pderiv{i}u(x)}
      + \sum_{1 \le i \le j \le d}
         \abs{\diff^2_{h,ij} u_h(x) - \pderiv[2]{ij} u(x)}\\
   &\le
      \sum_{i = 1}^d \abs{\diff_{h,i} (u_h(x) - u(x))}
         + \sum_{i = 1}^d \abs{\diff_{h,i} u(x) - \pderiv{i}u(x)}\\
   &\quad
      + \sum_{1 \le i \le j \le d} \abs{\diff^2_{h,ij} (u_h(x) - u(x))}
      + \sum_{1 \le i \le j \le d}
         \abs{\diff^2_{h,ij} u(x) - \pderiv[2]{ij} u(x)}\\
   &\le 2 d \varepsilon_1 + d \varepsilon_1 + \frac{\varepsilon_1}{c_d}.
\end{align*}

Thus, for some constant $c'_d$ depending only on $d$, the inequality~\eqref{equ:lemma:approx:3} holds with $\varepsilon$ replaced by $c'_d\,\varepsilon_1$.

The result follows now by taking $\varepsilon_1$ and $h_0$ appropriately.
\endproof}

The main implication of theorem~\ref{theorem:discrete:approx} is that any smooth convex function is a limit of a sequence of functions with positive semidefinite discrete Hessians, giving an affirmative answer to the first issue mentioned at the beginning of this section.
Moreover, an application of lemma~\ref{lemma:cont:1} implies this result also for non-smooth convex functions, with convergence in $\|\cdot\|_0$ on compact subsets of $\Omega$ (see definitions below).

We would like to show next that our set of approximants also solves the second issue.
That is, if we have a convergent (in certain norm) sequence of functions with positive semidefinite discrete Hessians, then the limit is convex.
On the other hand, if the sequence is not convergent, we would like also to understand under which further conditions on the sequence we may extract a subsequence converging to a convex function.

In what follows, we will work with sequences of functions defined on finer and finer meshes, that is, sequences $S = {(u_{h_n})}_{n\in\NN}$ with $u_{h_n}\in\U_{h_n}$ for every $n$, such that $h_n/h_{n+1}\in\NN$ and ${(h_n)}_n$ decreases to $0$.
We will denote by $\Seq$ the set of all such sequences, and use the notation
\[
   \mesh(S) = \cup_{n}\, \mesh_{h_n},
   \quad
   \interior{\mesh}(S) = \mesh(S)\cap\interior{\Q},
   \quad\text{and}\quad
   \bdy\mesh(S) = \mesh(S)\cap\bdy\Q.
\]
So as not to clutter even more the notation, usually we will drop the index $n$, writing $S = (u_h)$ for $S\in\Seq$, when this does not lead to confusion.

If $S = (u_h)\in\Seq$ and $u:\Q\to\RR$, we will say that
\[
   \lim_{h\to 0} u_h(x) = u(x)
   \quad\text{uniformly for $x\in\mesh(S)$},
\]
if for any $\varepsilon > 0$, there exists $h_0 = h_0(\varepsilon) > 0$, so that
\( \abs{u_h(x) - u(x)} < \varepsilon \)
for all $h < h_0$ and $x\in \mesh_{h}$, that is,
$\max_{x\in\mesh_h} \abs {u_h(x) - u(x)} \to 0$ as $h\to 0$.
In this case we will write, with a little abuse of notation,
\[ \lim_{h\to 0}\,\normast{u_h - u}_{0} = 0. \]

If in addition $u\in C^{1}$, we write, similarly,
\[ \lim_{h\to 0}\,\normast{u_h - u}_{1} = 0, \]
to indicate that
\( \lim_{h\to 0}\,\normast{u_h - u}_{0} = 0 \)
and
\[ \lim_{h\to0}\,\diff_{h,i} u_h(x) = \pderiv{i} u(x) \]
uniformly for all $x\in\mesh(S)$ for which the finite differences make sense and $i = 1,\dots,d$.

{\lemma\label{lemma:discrete:convergence}
Suppose $S = (u_{h}) \in\Seq$ and $u\in C^{1}$ are such that
\[ \lim_{h\to 0}\,\normast{u_{h} - u}_1 = 0, \]
and
\[ H_h u_h (x) \succeq 0 \quad\text{for all $x\in\interior{\mesh}(S)$}. \]
Then $u$ is convex.
\endlemma}

{\proof
We follow closely what was done in section~\ref{sec:convex:cont}, using a variant of the divergence theorem for discrete variables, so that only approximations to the function or its derivatives---but not the second derivatives---are needed.

For $x\in\interior{\mesh}_{h'}$ and $h\ll h'$, let us define
\[ H^{*}_{h',h}\, u_{h} (x) = \sum_{y} H_{h}\, u_h(y), \]
where the sum is over all $y\in\mesh_h$ such that $\Q_{h}(y)\subset \Q_{h'}(x)$.

Since the sum of positive semidefinite matrices is a positive semidefinite matrix, we have
\[
   H^{*}_{h',h}\, u_{h} (x) \succeq 0
   \quad\text{for all $x\in\interior{\mesh}_{h'}$.}
\]

In one dimension we have just one entry in $H^{*}_{h',h}$, involving a term of the form
\begin{multline*}
   u(x - h') - 2\,u(x - h' + h) + u(x - h' + 2h)\\
      + u(x - h' + h) - 2\, u(x - h' + 2h) + u(x - h' + 3h)
      +\dotsb\\
   + u(x + h' - 2h) - 2\, u(x + h' - h) + u(x + h')\\
   = \bigl(u(x + h') - u(x + h' - h)\bigr)
      - \bigl(u(x - h' + h) - u(x - h')\bigr).
\end{multline*}
That is, if $d = 1$ then
\begin{equation}\label{equ:pf:2}
   H^{*}_{h',h}\, u_{h} (x) =
      \begin{bmatrix}
      \frac{1}{h}\,
         \bigl(
            \diff_{h} u_{h} (x + h' - h) - \diff_{h} u_{h} (x - h')
         \bigr)
      \end{bmatrix}.
\end{equation}

If $d = 2$ and $x = (x_1,x_2)$, the diagonal entries are similar to the one dimensional case.
For instance,
\begin{multline}\label{equ:pf:3}
   \bigl(H^{*}_{h',h}\, u_{h} (x)\bigr)_{11} =
      \frac{1}{h}\,\Bigl(
      \sum_{k} \diff_{h,1} u_{h} (x_1 + h' - h, x_2 + kh)\\
      -
      \sum_{k} \diff_{h,1} u_{h} (x_1 - h', x_2 + kh)
      \Bigr),
\end{multline}
where both sums are on $k$'s such that $-h' < kh < h'$.
On the other hand, the off-diagonal terms are of the form
\begin{multline}\label{equ:pf:4}
   \bigl(H^{*}_{h',h}\, u_{h} (x)\bigr)_{12} =
      \frac{1}{4{h}^2}\,
      \Bigl(
         \bigl(
              u_{h} (x + \alpha)
            + u_{h} (x + \beta)
            + u_{h} (x + \gamma)
            + u_{h} (x + \delta)
         \bigr)\\
         -
         \bigl(
              u_{h} (x + \alpha')
            + u_{h} (x + \beta')
            + u_{h} (x + \gamma')
            + u_{h} (x + \delta')
         \bigr)\\
         -
         \bigl(
              u_{h} (x - \alpha')
            + u_{h} (x - \beta')
            + u_{h} (x - \gamma')
            + u_{h} (x - \delta')
         \bigr)\\
         +
         \bigl(
              u_{h} (x - \alpha)
            + u_{h} (x - \beta)
            + u_{h} (x - \gamma)
            + u_{h} (x - \delta)
         \bigr)
      \Bigr),
\end{multline}
where
\begin{align*}
   \alpha  &= (h' - h)\, \e_1 + (h' - h)\, \e_2, &
   \alpha' &= (h' - h)\, \e_1 - (h' - h)\, \e_2,\\
   \beta   &= h'\, \e_1 + (h' - h)\, \e_2, &
   \beta'  &= h'\, \e_1 - (h - h')\, \e_2,\\
   \gamma  &= (h' - h)\, \e_1 + h'\, \e_2, &
   \gamma' &= (h' - h)\, \e_1 - h'\, \e_2,\\
   \delta  &= h'\, \e_1 + h'\, \e_2, &
   \delta' &= h'\, \e_1 - h'\, \e_2.
\end{align*}

For $d > 2$ we get similar expressions to those of equations~\eqref{equ:pf:2} and~\eqref{equ:pf:3} for the diagonal terms, and~\eqref{equ:pf:4} for the off-diagonal terms, except that---as when going from~\eqref{equ:pf:2} to~\eqref{equ:pf:3}---they must be summed over mesh points on $d-1$ dimensional surfaces perpendicular to the (say) $i$-th direction, or $d-2$ dimensional surfaces perpendicular to the (say) $ij$ plane.

This was also the case in the equations~\eqref{equs:cont:2} for the continuous case, and so we can think of these sums as approximations to those integrals, where small $d$-dimensional cubes of side length $h$ centered at mesh points have been used.
In order to do this we must multiply the entries by $h^{d}$ so as to obtain the correct dimensions.

In fact, since we are assuming
\[ \lim_{h\to 0} \normast{u_h - u}_1 = 0, \]
and $H^{*}_{h',h} u_h(x)$ involves either no finite differences (off the diagonal) or only first order differences (in the diagonal) which converge uniformly to (respectively) $u$ or its first derivatives,
\[
   \lim_{h\to 0}\, h^{d}\,H^{*}_{h',h} u_h (x)
   = R_{h'} u (x)
   \quad\text{for all $x\in\interior{\mesh}_{h'}$,}
\]
where
\[ R_{h'} u (x) \succeq 0, \]
since the limit of positive semidefinite matrices is positive semidefinite.

On the other hand, due to the convergence of $u_h$ to $u$ in the $\|\cdot\|_1^*$ norm, we can see that, as $h\downarrow 0$,
\begin{align*}
\bigl(h^{d}\,H^{*}_{h',h} u_h (x) \bigr)_{\!ii}
&\to
       \int_{F^{+}_{\delta,i}(x)} \pderiv{i} u(y)\,\dd{y}
     - \int_{F^{-}_{\delta,i}(x)} \pderiv{i} u(y)\,\dd{y},
\\
\bigl(h^{d}\,H^{*}_{h',h} u_h (x)\bigr)_{\!ij}
&\to
   \Bigl(\int_{F^{+}_{\delta,i}(x)\cap F^{+}_{\delta,j}(x)} u(y)\,\dd{y}
      - \int_{F^{+}_{\delta,i}(x)\cap F^{-}_{\delta,j}(x)} u(y)\,\dd{y}
   \\[4pt]
&\qquad    -
    \int_{F^{-}_{\delta,i}(x)\cap F^{+}_{\delta,j}(x)} u(y)\,\dd{y}
      + \int_{F^{-}_{\delta,i}(x)\cap F^{-}_{\delta,j}(x)} u(y)\,\dd{y}
      \Bigr),
\end{align*}
if $i\neq j$.
Thus,
\[
   R_{h'} u(x) = (2h')^d\, \Hw_{h'} u(x)
   \quad\text{for $x\in\interior{\mesh}_{h'}$},
\]
which implies
\( \Hw_{h'} u(x) \succeq 0 \)
for all $x\in\interior{\mesh}_{h'}$.
The result now follows from theorem~\ref{theorem:cont:1}.
\endproof}

The convergence in $\normast{\cdot}_1$ should be guaranteed by the definition of the discrete problem and is problem dependent.
We make now a general assumption on an algorithm to ensure convergence in $\normast{\cdot}_1$.

Let us denote by \(\Lambda^1\) the space of Lipschitz continuous functions defined on $\Q$, recalling that \(\Lambda^1 = W^{1,\infty}\), the space of continuous functions having derivatives in the weak sense up to order one bounded.
If $K > 0$, we let $\Lambda^1_K$ be the set of functions $u\in \Lambda^1$ such that both $\abs{u(x)} \le K$ for all $x\in\Q$ and \(\abs{u(x) - u(x+t\,\e_i)}\le Kt\) for all $i = 1,\dots,d$ and $t > 0$ whenever $x,x + t\,\e_i\in\Q$.
Of course, the condition \(\abs{u(x) - u(x+t\,\e_i)}\le Kt\) is equivalent to $\abs{\pderiv{i}u(x)}\le K$.

Similarly, let us denote by \(\Lambda^2 = W^{2,\infty}\) the space of functions whose first derivatives are in $\Lambda^1$, and by $\Lambda^2_K$ the set of all functions in $\Lambda^2$ which have all weak derivatives up to order $2$ bounded by $K$.
In particular, we have $\Lambda^2\subset C^1$.

By analogy to the continuous case, let us define for $h$ and $K$ positive the following spaces of discrete functions:
\begin{equation}\label{equ:wh:defs}
\begin{gathered}
   \Lambda^{0}_{h,K}
      = \{u \in\U_h : \abs{u(x)} \le K, x\in\mesh_h \},\\
   \Lambda^{1}_{h,K}
      = \{u \in \Lambda^{0}_{h,K} :
         \abs{\diff_{h,i} u(x)} \le K, x\in\mesh_h, i = 1,\dots,d \},\\
   \Lambda^{2}_{h,K}
      = \{u \in \Lambda^{1}_{h,K} :
         \abs{\diff^2_{h,ij} u(x)} \le K,
         x\in\mesh_h,
         i,j = 1,\dots,d \},
\end{gathered}
\end{equation}
with the understanding that on the right hand sides we take $\diff_{h,i}u(x)$ or $\diff^2_{h,ij}$ for all $x,i$ and $j$ where it makes sense.

The following result is a version of the Arzela-Ascoli theorem, and we omit its proof.

{\theorem\label{theorem:arzela}
If $S = (u_h)\in\Seq$ for which $u_h\in \Lambda^{1}_{h,K}$  for all $h > 0$, then there exists $u\in \Lambda^{1}_{K}$ and a subsequence $S' = (u_{h'})$ of $S$ such that
\[ \lim_{h'\to 0}\, \normast{u_{h'} - u}_0 = 0. \]
\endtheorem}

Combining the previous results, we have:

{\theorem\label{theorem:discrete:converge}
If $S = (u_h)\in\Seq$ is such that
\[
   u_h\in \Lambda^{2}_{h,K}
   \quad\text{and}\quad
   H_h u_h (x)\succeq 0
\]
for all $x\in\interior{\mesh}(S)$ and $h > 0$, then there exists $u\in \Lambda^{2}_K$, and a subsequence $S' = (u_{h'})$ of $S$ such that
\[
   \lim_{h'\to 0}\, \normast{u_{h'} - u}_1 = 0
   \quad\text{and}\quad
   \text{$u$ is convex}.
   \]
\endtheorem}

{\proof
Applying theorem~\ref{theorem:arzela} to the functions $\diff_{h,i}$, perhaps on some smaller meshes $\mesh'_{i,h}$ with $\interior{\mesh}_h \subset \mesh'_{i,h} \subset \mesh_h$, we may find functions $u_i\in \Lambda^{1}_K$ and a subsequence of $S$, $S' = (u_{h'})$, such that
\[
   \lim_{h'\to 0}\,\normast{\diff_{h',i} u_{h'} - u_i}_0 = 0
   \quad\text{for $i = 1,\dots,d$.}
\]

To show that there exists $u$ such that $u_i = \pderiv{i} u$ (in the classical sense), we define for $x = (x_1,\dots,x_d)\in\Q$,
\[
   u\Bigl(\sum_{i = 1}^d x_i\,\e_i\Bigr)
   = u(\vvv{0})
      + \sum_{j = 1}^d
         \int_{0}^{x_j}
            u_j\Bigl(t_j\, \e_j + \sum_{i = 1}^{j-1} x_i\,\e_i\Bigr)
         \,\dd{t_j}.
\]

Since for $i = 1,\dots,d$, $u_i$ is bounded (by $K$) and continuous, $u$ is well defined and continuous.
Using that $\diff_{h',i} u_{h'}$ converge uniformly to $u_i$, for $x\in\mesh(S')$ we may write
\[
   u(x) =
   \lim_{h'\to 0}\,
      u_{h'}(\vvv{0})
      + \sum_{j = 1}^d
         \sum_{k_j = 0}^{x_j/h'} h'
            \diff_{h',j}
            u_{h'}\Bigl(k_j h'\, \e_j + \sum_{i = 1}^{j-1} x_i\,\e_i\Bigr).
\]

But, for a given $h'$ and all $x\in\mesh_{h'}$ we have
\begin{multline}\label{equ:approx:indep}
   u_{h'}(\vvv{0})
   + \sum_{j = 1}^d
      \sum_{k_j = 0}^{x_j/h'} h'
         \diff_{h',j}
         u_{h'}\Bigl(k_j h'\, \e_j + \sum_{i = 1}^{j-1} x_i\,\e_i\Bigr)
   \\=
   u_{h'}(\vvv{0})
   + \sum_{j = 1}^d
      \biggl(
         u_{h'}\Bigl(\sum_{i = 1}^{j} x_i\,\e_i\Bigr)
         -
         u_{h'}\Bigl(\sum_{i = 1}^{j-1} x_i\,\e_i\Bigr)
      \biggr)
   =
      u_{h'}(x),
\end{multline}
and therefore
\[ \lim_{h'\to 0} \normast{u_{h'} - u}_0 = 0. \]

Arguing as in equation~\eqref{equ:approx:indep} and taking limits, we may also verify the ``independence of path'', that is, for $x\in\mesh_{h'_0}$ and $h' \le h'_0$ we may write,
\[
   u(x + h'\,\e_i) - u(x) = \int_{0}^{h'} u_i(x + t_i\,\e_i)\,\dd{t_i}
   \quad\text{for all $i = 1,\dots,d$.}
\]

Using the continuity of $u$, we see that the last equation is valid for all $x\in\interior{\Q}$ and $h'$ small enough, and, using once more the continuity of $u$, that
\[
   u(x + \delta\,\e_i) - u(x)
   = \int_{0}^{\delta} u_i(x + t_i\,\e_i)\,\dd{t_i}
   \quad\text{for all $i = 1,\dots,d$,}
\]
for all $x\in\interior{\Q}$ and all $\delta > 0$, and therefore
\[
   \pderiv{i} u(x) = u_i(x)
   \quad\text{for all $x\in\interior{\Q}$ and $i = 1,\dots,d$.}
\]

The convexity of $u$ follows from lemma~\ref{lemma:discrete:convergence}.
\endproof}

{\remark\label{remark:noconvergenceof2}
If $u_n$ is a sequence of convex functions converging pointwise to $u$, then $u$ must be convex.
However it is not true in general that the second derivatives or the Hessian of $u_n$ will converge to that of $u$, even if we have uniform convergence of $u_n$ and its derivatives to those of $u$.
In this sense, theorem~\ref{theorem:discrete:converge} cannot be bettered too much.

For example, consider in $d = 1$ the functions
\[ u(x) = \frac{x^2}{2}, \]
and, with $h = 1/n$ for $n\in\NN$,
\[ u_h(x) = u(x) + \frac{\cos (\pi n x)}{{(n\pi)}^2}. \]

We have,
\[
   \begin{array}{cc}
   u_h(x) \ge 0,
   &   \lim_{h\to 0}\,u_h(x) = u(x),\\[4pt]
   0 \le u_h'(x) = u'(x) - \dfrac{\sin (\pi n x)}{n\pi} \le 1,
   &   \lim_{h\to 0}\,u'_h(x) = u'(x),
   \end{array}
\]
where the limits are uniform in $x$.
Also,
\begin{gather*}
   u_h''(x) = u''(x) - \cos(\pi n x)\ (\ge 0),
   \qquad\text{and}\qquad
   \Delta^2 u_h (x) = 1 - \frac{4 \cos (\pi n x)}{\pi ^2}.
\end{gather*}

Taking $n = 2^{j + m}$, $x = \frac{i}{2^j} = \frac{2^m\,i}{n}$, for $m > 0$ we obtain
\[
   \Delta^2 u_h (x)
   = 1 - \frac{4 \cos (\pi 2^m i)}{\pi ^2}
   = 1 - \frac{4}{\pi^2} \approx 0.594715,
\]
and so the second order differences at dyadic points converge to this constant value.
However, also at these points,
\( u''(x) = 1 \).\END
\endremark}

\smallskip
On the other hand, we may weaken the conditions on convergence, following what was done for the one-dimensional case in section~\ref{sec:review}.

If $u_h$ is defined in $\mesh_h$ and $H_h u_h(x)\succeq 0$ for $x\in\mesh_h$, its restriction to the one-dimensional line obtained by fixing all coordinates except the $i$-th one, satisfies
\[
   \frac{u(x + i h\, \e_i) - u(x + (i-1)\, h\, \e_i)}{h}
   \le \frac{u(x + (i+1)\, h\,\e_i) - u(x + i h\, \e_i)}{h},
\]
for $i = 1,\dots,d$, since these are coefficients in the main diagonal of $H_h u_h(x + h\e_i)$.
That is, the inequalities~\eqref{equ:1d:convex} are satisfied, and hence if for some $M > 0$,
\[ \abs{u_h(x)} \le M \quad\text{for all $x\in\mesh_h$,} \]
then for any given $\varepsilon > 0$ we may find a constant $C$, depending on $M$ and $\varepsilon$ but not on $u$, such that
\[
   \abs{\diff_{h,i} u(x)} \le C
   \quad\text{for all $x\in \mesh_h$ such that $\dist(x,\bdy\Q) > \varepsilon$.}
\]

Therefore, by applying theorem~\ref{theorem:arzela} to $d$-dimensional cubes contained in $\interior{\Q}$, we may strengthen lemma~\ref{lemma:discrete:convergence} to obtain a variant of lemma~\ref{lemma:1d:convergence:1}:

{\corollary\label{corollary:discrete:convergence}
Let $S = (u_{h})\in\Seq$ such that
\[ H_h u_h (x) \succeq 0 \quad\text{for all $x\in\interior{\mesh}(S)$}, \]
and suppose $u$ is a (finite) function defined on $\mesh(S)$ satisfying
\[
   \lim_{h\downarrow 0} u_h(x) = u(x)
   \quad\text{for all $x\in \mesh(S)$.}
\]
Then, the convergence of $u_h$ to $u$ is uniform on compact subsets of $\interior{\Q}$, and $u$ may be uniquely extended to a convex function defined on all of $\Q$.
\endcorollary}

Similarly, we may weaken the conditions of theorem~\ref{theorem:arzela} to obtain the following version of lemma~\ref{lemma:1d:convergence:2}:

{\corollary\label{corollary:arzela}
If $M > 0$ and $S = (u_h)\in\Seq$ are such that
\[ \abs{u_h(x)} \le M \quad\text{for all $h$ and $x\in \mesh_h$}. \]
Then there exists a continuous convex function $u$ defined on $\Q$ and a subsequence $S' = (u_{h'})$ of $S$ such that
\[
   \lim_{h'\to 0}\, u_{h'}(x) = u_{h}(x)
   \quad\text{for all $x\in\mesh(S')$}.
\]
Moreover, the convergence is uniform on compact subsets of $\interior{\Q}$.
\endcorollary}

We conclude this section with some comments.

{\remark[boundary behavior]\label{remark:bdy}
If $u\in C^2$, $\pderiv[2]{ij} u(x)$ is defined initially for $x\in\interior{\Q}$, but the very definition of $C^2$ as the set of those functions in $C^2(\interior{\Q})$ having second order derivatives uniformly continuous on $\interior{\Q}$, makes it possible to define $\pderiv[2]{ij} u(x)$ for $x\in\bdy\Q$ by a limiting argument.
And the same goes for lower order derivatives, and even $Hu(x)$.

The situation is different for discrete functions defined on $\mesh_h$ for fixed $h > 0$, since we cannot take limits.
However, the particular geometry of $\Q$ makes it possible (as we did) to consider for $x\in\bdy\mesh_h$ as many finite differences as we can.
For instance, $\diff^2_{ii} u(x)$ is well defined for $x\in\bdy\mesh_h$ as long as $x \pm h\,\e_i\in\mesh_h$.

Pushing the definitions a little further, we may define the discrete Hessian $H_h u(x)$ for $x\in\bdy{\mesh}_h$, by including as many second derivatives as we can.
For example, if $d = 3$ and $h = 1/2$, we can define
\[
   \begin{array}{cll}
   H_h u(x)
   =
      \begin{bmatrix}
      \diff^2_{h,11} u(x) & \diff^2_{h,12} u(x)\\[2pt]
      \diff^2_{h,12} u(x) & \diff^2_{h,22} u(x)\\[2pt]
      \end{bmatrix}
   &\quad
   &\text{if $x = (1/2,1/2,0)$},\\[14pt]
   H_h u(x)
   =
      \begin{bmatrix}
      \diff^2_{h,11} u(x)
      \end{bmatrix}
   &&\text{if $x = (1/2,0,0)$},
   \end{array}
\]
leaving $H_h u(x)$ undefined if $x = (0,0,0)$.

As the reader may verify, our previous results involving $H_h$ remain valid with this interpretation of the discrete Hessian.\END
\endremark}

\section{Approximating functionals}
\label{sec:functionals}

We are in position now to use finite difference approximations of a wide class of optimization problems on convex functions.

Let us describe this technique by assuming, for instance, that $(V,\norm{\cdot}_V)$ is a Banach space of real valued functions on $\Q$, the functional
\[ J(v) = \int_\Q F(x,v(x),\grad v(x)) \,\dd{x} \]
is defined and continuous on $V$, and we are interested in the optimization problem
\begin{equation}\label{prob:PJ}
   \inf\, \{ J(v): v\in \Conv\},
\end{equation}
where $\Conv$ is a family of convex functions, $\Conv\subset V$.

If the functions in $\Conv$ may be approximated by convex functions in $C^3\cap V$, then using theorem~\ref{theorem:discrete:approx} it may be not too difficult to define for each $h > 0$ (or a sequence converging to $0$ of such $h$'s), a family $\Conv_h$, $\Conv_h\subset\U_h$, and a functional $J_h$ defined on $\Conv_h$, such that:
{\enumerate
\item $H_h v_h(x) \succeq 0$ for all $v_h\in\Conv_h$ and $x\in\interior{\mesh}_h$,

\item for any $v\in \Conv$ and any $\varepsilon > 0$, there exists $h > 0$ and $v_h\in \Conv_h$ such that $\abs{J_h(v_h) - J(v)} < \varepsilon$.
\endenumerate}

Condition 2 immediately implies that
\begin{equation}\label{equ:fund:1}
   \inf\, \{ J(v): v\in \Conv\} \ge \inf_{h}\,\inf\, \{ J_h(v_h) : v_h\in \Conv_h \}.
\end{equation}

To prove the converse, we observe that
\[
   \inf_h\, \{J_h(v_h) : v_h\in \Conv_h \}
   = \inf_{h,K > 0}\, \inf\, \{J_h(v_h) : v_h\in \Conv_{h,K}\},
\]
where $\Conv_{h,K} = \Conv_h\cap \Lambda^{2}_{h,K}$.
Keeping $K > 0$ fixed, we may find a sequence $S = (u^K_h)$, with $u^K_h\in\Conv_{h,K}$ such that
\[
   J_h(u^K_h)\downarrow \inf_{h}\, \inf\,\{J_h(v_h) : v_h\in \Conv_{h,K}\},
\]
and using theorem~\ref{theorem:discrete:converge}, a subsequence $S' = (u^K_{h'})$ and a convex function $u^K\in C^{2}$ with $\normast{u^K_h - u^K}_1$ converging to $0$.

If $V_h$ and $J_h$ are such that
\[
   u^K\in \Conv
   \quad\text{and}\quad
   J_h(u^K_h)\to J(u^K),
\]
letting $K \to \infty$ we will have
\begin{equation}\label{equ:fund:2}
   \inf\,
   \{ J(v): v\in \Conv\} \le \inf_{h}\,\inf\,\{ J_h(v_h) : v_h\in \Conv_h \}.
\end{equation}

Putting together the inequalities~\eqref{equ:fund:1} and~\eqref{equ:fund:2}, we will have discrete approximations of the problem~\eqref{prob:PJ}.

Let us give some more concrete examples.
Suppose, for instance that $V = H^1$ is the set of functions $u: \Q\to\RR$ with finite norm
\[
   \norm{u}
   = \Bigl(\int_\Q
      {\abs{\grad u (x)}}^2 + {\abs{u(x)}}^2\,\dd{x}
      \Bigr)^{1/2},
\]
and suppose $\Conv$ is the set of all convex functions in $H^1$.
Given $f\in V$ we would like to find its projection on $\Conv$, that is, find $u\in\Conv$ such that
\[ \norm{u - f} = \min_{v\in\Conv} \norm{v - f}, \]
and, getting rid of the square roots, we may set
\[ J(v) = {\norm{v - f}}^2. \]

In this example we actually have a unique minimum, since the norm is strictly convex.
We may consider then $\Conv_h$ as the set of discrete functions $v_h\in\U_h$ with $H_h v_h(x)\succeq 0$ for all $x\in\interior{\mesh}_h$.

Assuming for simplicity that $f\in C^1$, for $v_h\in\U_h$ we define
\[
   J_h(v_h)
   = h^d\,\Bigl(
      \sum_{x\in\interior{\mesh}_h} \bigl(
         {\abs{v_h (x) - f(x)}}^2
         +
         \sum_{i=1}^d {\abs{\diff_{h,i} v_h (x) - \pderiv{i} f(x)}}^2
         \bigr)
      \Bigr).
\]

Then, it is easy to see that~\eqref{equ:fund:1} and~\eqref{equ:fund:2} hold.
In fact, the convex functions of $V$ may be approximated by $C^3$ convex functions, so that, given $\varepsilon > 0$ there exists $v\in C^3\cap H^1$ such that
\[ \norm{u - v} < \frac{1}{2}\,\varepsilon, \]
and use theorem~\ref{theorem:discrete:approx} to find $h$ small enough and $v_h\in\Conv_h$ such that
\[ \abs{J_h(v_h) - J(v)} < \frac{1}{2}\,\varepsilon. \]
Thus we will have
\[
   \abs{J_h(v_h) - J(u)}
   \le \abs{J_h(v_h) - J(v)} + \abs{J(v) - J(u)}
   < \varepsilon,
\]
for some $v_h\in\Conv_h$.

\section{Numerical results}
\label{sec:numerical}

In this section we illustrate the behavior of the numerical scheme by applying it to the problems mentioned in the introduction, namely,
the monopolist problem,
norm projections on the set of convex functions,
and fitness of data by discrete convex functions.

In all of the following examples, we associate with each mesh node $P_k$, $k = 1,\dots,N$, the unknown value $u_k$, and the square $Q_k = \Q\cap \Q_h(P_k)$, of area $\abs{Q_k}$.
We also consider the discrete Hessian $H_h(P_k)$ as discussed in the remark~\ref{remark:bdy}, i.e., imposing convex constraints on the boundary whenever they make sense.

Even though theorem~\ref{theorem:discrete:converge} requires us to impose upper bounds on the second order differences in order to ensure convergence, the following numerical experiments were carried on without this requirement.
Convergence was nevertheless observed at optimal rates.

The times reported correspond to the experiments being run on a PC, with a 2.8GHz Pentium~IV processor and 2GB of RAM, running \emph{Linux}.
The matrices were assembled using \emph{OCTAVE}~\cite{octave} and the semidefinite program was solved using \emph{CSDP}~5~\cite{csdp1} with the default parameters.
The graphics were obtained using \emph{Mathematica}~\cite{mathematica}.

\subsection{The monopolist problem}\label{sec:economics}

Since this problem is not widely known in the mathematics community, let us start by giving a brief description of it following the one given in~\cite{MV07}, where it is referred to as the \emph{revenue maximization in a multiple-good monopoly} problem.

{\problem[The monopolist problem]\label{prob:monopolist}
A seller with $d$ different objects faces a single buyer, whose preferences over consumption and money transfers are given by $U(x,p,t) = x\cdot p - t$, where $x\in [0,1]^d$ is the vector of the buyer's valuations, $p\in\RR^d$ is the vector of quantities consumed for each good, and $t\in\RR$ is the monetary transfer made to the seller.
The valuations $x$ are only observed by the buyer, and a density function $f(x)$ represents the seller's belief on the buyer's private information $x$.
The seller's problem is to design a revenue maximizing mechanism to carry out the sale, and it is enough to consider only direct revelation mechanisms: the buyer must prefer to reveal its information truthfully (incentive compatibility) and to participate voluntarily (individual rationality).
Under these conditions, it can be proved that the seller's problem may be written as
\[
   \max_{u\in\Conv}
      \int_\Q \bigl( \grad u(x) \cdot x - u(x) \bigr) \,f(x)\, \dd{x},
\]
where
$\Q = [0,1]^d$,
$f$ is a non-negative probability density function over $\Q$,
and $\Conv$ is the set of functions $u$ satisfying
{\enumerate
\item $u$ is convex,
\item $0 \le \grad u(x) \le 1$ for all $x\in \Q$ (the gradient taken in the weak sense and the inequalities componentwise), and
\item $u(0) = 0$.
\endenumerate}

The functional to be maximized is the seller's expected revenue.
For a buyer of type $x$, the solution $u$ to this optimization problem represents the utility received by her, and the $i$-th  component of $\grad u$ denotes the probability that she will obtain good $i$.
The restriction of convexity stems from incentive compatibility, and the condition $u\ge 0$ from individual rationality.\END
\endproblem}

We show now some numerical results for the problem~\ref{prob:monopolist} for the special case $f = 1$, for which analytic solutions in 2 and 3 dimensions are known, allowing us to judge the behavior of the discrete approximations.
That is, the functional to be minimized is
\[ J(u) = \int_Q \bigl( u(x) - \grad u(x) \cdot x \bigr)\, \dd{x}. \]

In both 2 and 3 dimensions, the solutions are piecewise linear convex functions whose partial derivatives are either $0$ or $1$.
For example, for $d = 2$ the solution is
\[ u(x_1,x_2) = \max\,\{0, x_1 - a, x_2 - a, x_1 + x_2 - b\}, \]
where
\[
   a = \frac{2}{3}
  \qquad\text{and}\qquad
  b = \frac{1}{3}\,\bigl(4-\sqrt{2}\bigr),
\]
and the value at the optimum is
\[ J(u) = \frac{2}{27}\, \bigl(6 + \sqrt{2}\bigr) \approx 0.549201. \]

\begin{figure}\centering
\begin{minipage}[b]{.48\textwidth}\centering
\includegraphics[width=.95\textwidth]{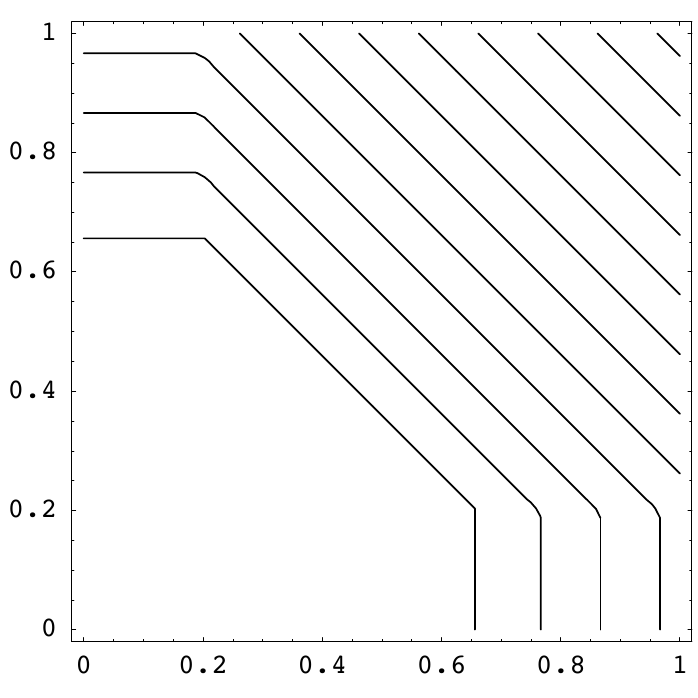}\\
(a)
\end{minipage}\hfill
\begin{minipage}[b]{.48\textwidth}\centering
\includegraphics[width=.95\textwidth]{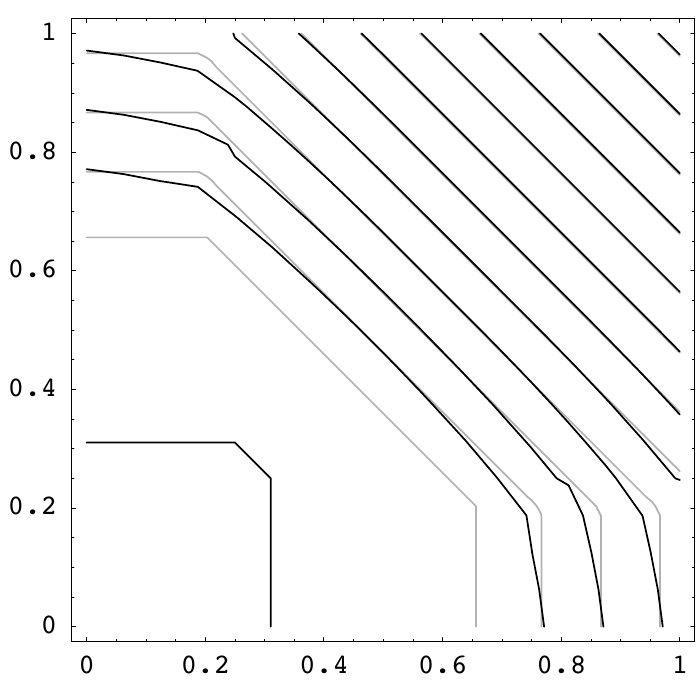}\\
(b)
\end{minipage}\\[10pt]
\begin{minipage}[b]{.48\textwidth}\centering
\includegraphics[width=.95\textwidth]{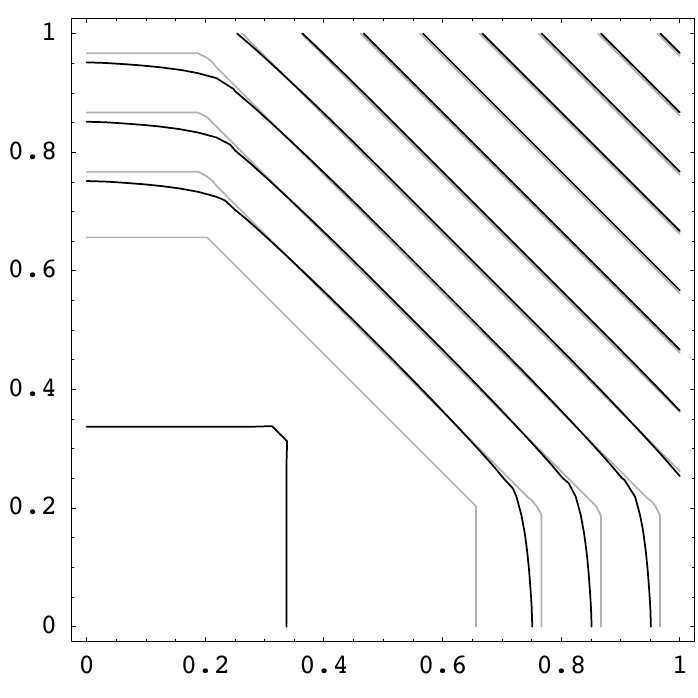}\\
(c)
\end{minipage}\hfill
\begin{minipage}[b]{.48\textwidth}\centering
\includegraphics[width=.95\textwidth]{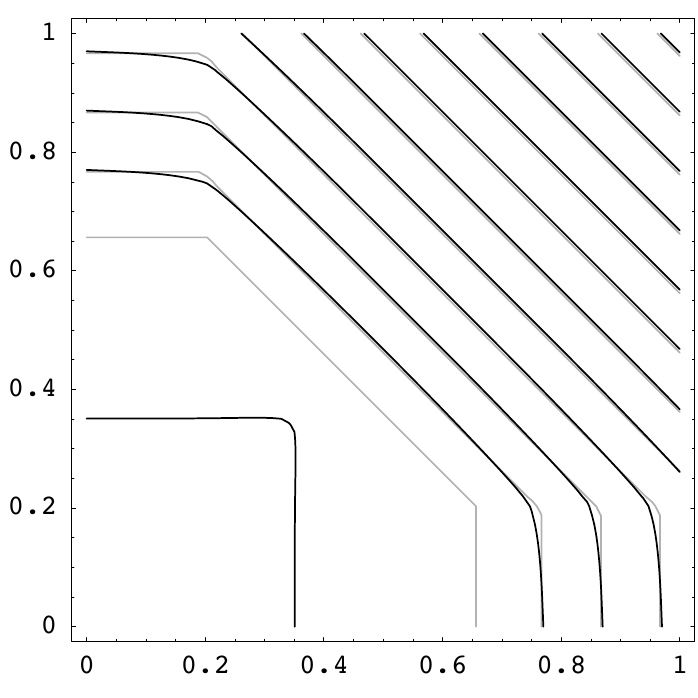}\\
(d)
\end{minipage}
\caption{%
Contour lines of the solutions of the economics problem in 2D:
analytic solution (a),
and discrete approximations with
$h = 1/16$ (b),
$h = 1/32$ (c),
$h = 1/64$ (d).
}
\label{fig:econo:2d}
\end{figure}

In figure~\ref{fig:econo:2d} we show the contour lines obtained.
In (a) the analytic solution, and then the discrete solution for $h = 1/16$, $1/32$ and $1/64$ in, respectively, (b), (c) and (d), with the contours of the analytic solution shown in a lighter gray.
The contours are $10^{-7},0.1,0.2,\dots,1.1$ (the CSDP solution is always positive due to the way it is solved).

We observe that the scheme introduces quite a bit of diffusion at lower resolutions, and, in particular, the heights at the point $(1,1)$ increase to the analytic solution at that point as $h$ decreases.
However, the jump in the gradients is well captured, especially at higher resolutions, even though we are asking for conditions on the discrete Hessians which are unbounded from above.

In table~\ref{tab:econo:2d} we give some comparative results at different mesh sizes.
In this table we indicate by $J_h(I_h u)$ the discrete functional evaluated at $I_h(u)$, the interpolant in $\mesh_h$ of the exact solution, and the error column refers to $\max\,\abs{u_h - I_h u}$ ($L^\infty$ error).

It is interesting to notice that even though the solution is not smooth ($u$ is only Lipschitz) the error in the $L^\infty$ norm is smaller or approximately equal to $h$.
In the last column of the table we show the quantity $(J(u)-J_h(u_h))/h$ which is approximately $0.14$ for all the values of $h$ reported, showing that even under such low regularity assumptions on $u$, the error in the functional behaves as $O(h)$.

\begin{table}\centering
\begin{tabular}{rlllcrc}
$n$ & $h = 1/n$ & $J_h(u_h)$ & $J_h(I_h u)$ & Error ($L^\infty$)
   & \multicolumn{1}{c}{Time}
   & $\frac{J(u)-J_h(u_h)}{h}$\\[2pt]
\hline\rule{0pt}{12pt}%
  8 & \ 0.125    & 0.5319 & 0.5444 & 0.0769 &   0.190s & 0.14 \\
 16 & \ 0.0625   & 0.5404 & 0.5478 & 0.0300 &   0.990s & 0.14 \\
 32 & \ 0.03125  & 0.5449 & 0.5488 & 0.0336 &  17.000s & 0.14 \\
 64 & \ 0.015625 & 0.5470 & 0.5491 & 0.0174 & 751.100s & 0.14 \\
 $\infty$ & \ 0 & 0.5492
\end{tabular}
\caption{Comparison of the approximations to the economics problems in 2D for different mesh sizes.
The column labeled $n$ denotes the number of subdivisions of the interval $[0,1]$ in each direction.
The number of unknowns is $(n+1)^2$.}
\label{tab:econo:2d}
\end{table}

In three dimensions there is also a solution with partial derivatives which are either $0$ or $1$, of the form
\begin{multline*}
   u(x,y,z) = \max\,
      \{0,
        x - a, y - a, z - a,\\
        x + y - b, x + z - b, y + z - b,
        x + y + z - c\}.
\end{multline*}

This time the coefficients cannot be expressed in a simple form, since they involve roots of polynomials of high degree, and we just give numerical approximations:
\[
   a = 0.840627,\quad
   b = 1.038352,\quad
   c = 1.236077,
\]
so that $b - a = c - b$, and the value of the functional is
\[ J(u) \approx 0.868405. \]

\begin{figure}\centering
\begin{minipage}[b]{.28\textwidth}\centering
\includegraphics[width=.94\textwidth]{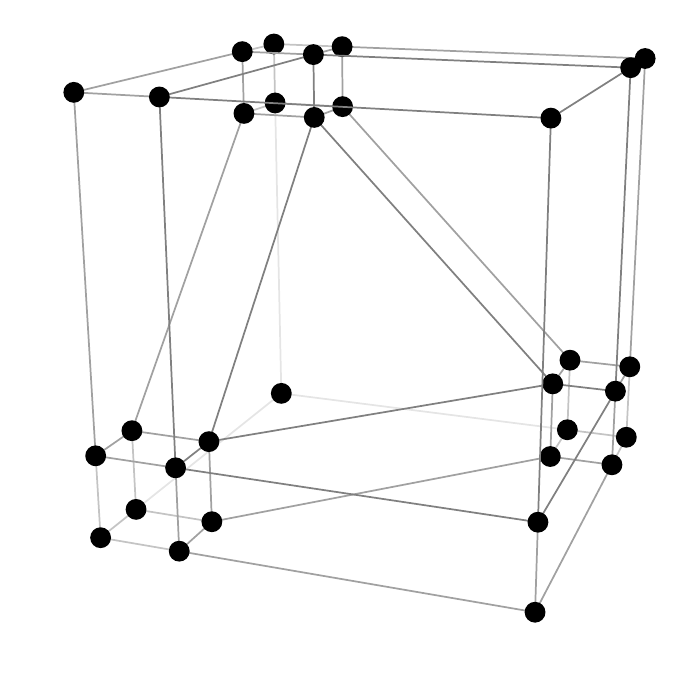}\\
(a)
\end{minipage}\hfill
\begin{minipage}[b]{.28\textwidth}\centering
\includegraphics[width=.94\textwidth]{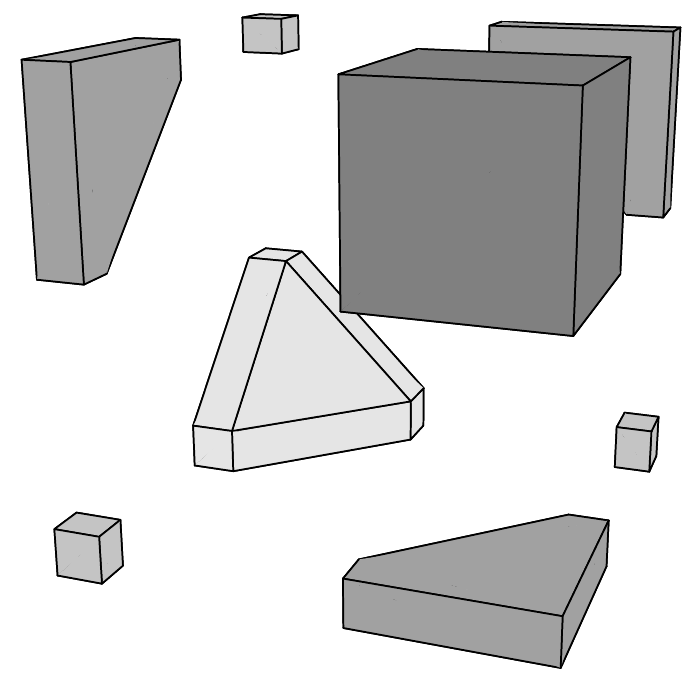}\\
(b)
\end{minipage}\hfill
\begin{minipage}[b]{.28\textwidth}\centering
\includegraphics[width=.94\textwidth]{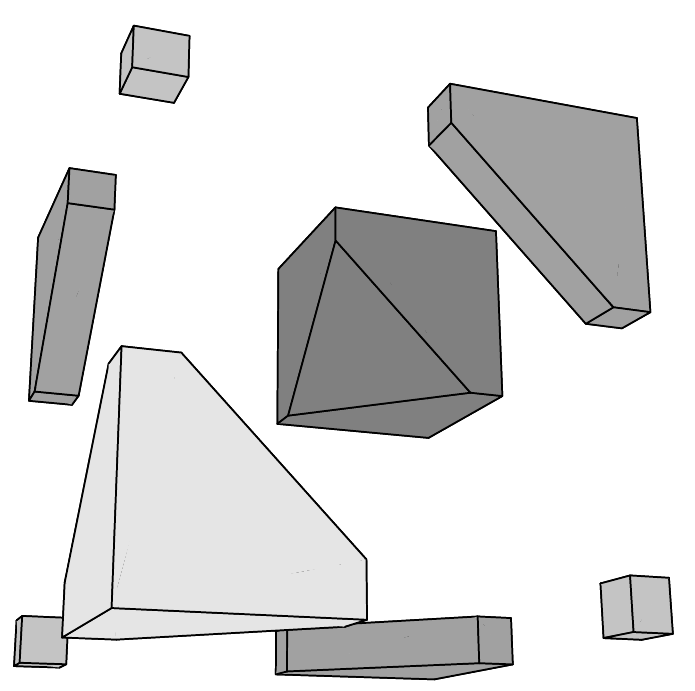}\\
(c)
\end{minipage}
\caption{%
Regions where the analytic solution in 3D has constant gradient:
with wire frames viewed from the positive octant (a),
and exploded views of the solid regions from the positive octant (b),
and the negative octant (c).
}
\label{fig:econo3D}
\end{figure}

In figure~\ref{fig:econo3D} we show the regions where $\grad u\cdot (1,1,1)$ is $0, 1, 2$ or $3$ in, different shades of gray.
We illustrate the regions with wire frames viewed from the positive octant in (a), and exploded views of the solid regions from the positive octant in (b) and the negative octant in (c).

In table~\ref{tab:econo:3d} we give some comparative results at different mesh sizes.
It is interesting to notice here that the $L^\infty$ error is not converging to zero with order $O(h)$.
Nevertheless, the  quantity $(J(u)-J_h(u_h))/h$ (shown in the last column) decreases slowly as $h$ goes to $0$, meaning that the error $J(u)-J_h(u_h)$ in the functional behaves as $O(h)$ in this range, exhibiting a similar behavior to that of the two dimensional case.

\begin{table}\centering
\begin{tabular}{rlcccrc}
$n$ & $h = 1/n$ & $J_h(u_h)$ & $J_h(I_h u)$ & Error ($L^\infty$)
   & \multicolumn{1}{c}{Time}
   & $\frac{J(u)-J_h(u_h)}{h}$\\[2pt]
\hline\rule{0pt}{12pt}%
  4 & \ 0.25     & 0.8195 & 0.8449 & 0.1356 &   0.21s & 0.20 \\
  8 & \ 0.125    & 0.8484 & 0.8605 & 0.1281 &   3.87s & 0.16 \\
 12 & \ 0.0833   & 0.8578 & 0.8647 & 0.1130 &  10.29s & 0.13 \\
 16 & \ 0.0625   & 0.8622 & 0.8661 & 0.1135 &  7m 20s & 0.10 \\
 20 & \ 0.0500   & 0.8648 & 0.8671 & 0.1177 & 50m 23s & 0.07 \\
 $\infty$ & \ 0 & 0.8684
\end{tabular}
\caption{%
Comparison of the approximations to the economics problems in 3D for different mesh sizes.
As before, the column labeled $n$ denotes the number of subdivisions of the interval $[0,1]$ in each direction.
The number of unknowns is $(n+1)^3$.}
\label{tab:econo:3d}
\end{table}

{\remark
Theoretically, semidefinite programs are polynomial time solvable.
However, as the tables~\ref{tab:econo:2d} and~\ref{tab:econo:3d} show, in practice we cannot go too far with the number of unknowns.

A reasonable size for the two dimensional problems we considered is a mesh of about $n = 40$ subdivisions.\END
\endremark}

\subsection{Projections}\label{sec:projections:2}

Carlier, Lachand-Robert and Maury~\cite{carlier} gave several examples of $H^1$ and $H_0^1$ projections, and in this section we consider one of the functions they considered, namely,
\[
   f(x_1,x_2)
   = -(4 + 5 x_1 x_2^2)\,e^{-30\,\left((x_1-1/2)^2 + (x_2-1/2)^2\right)}
   \quad\text{for $(x_1,x_2)\in\Q$.}
\]

\begin{figure}\centering
\includegraphics[width=.4\textwidth]{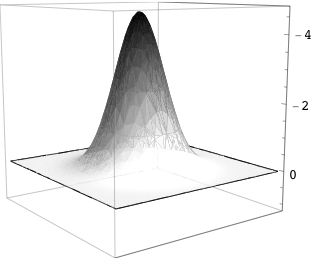}
\caption{Graph of $f(x_1,x_2) = -(4 + 5 x_1 x_2^2)\,e^{-30[(x_1-1/2)^2 + (x_2-1/2)^2]}$, shown upside down.}
\label{fig:fig2-f}
\end{figure}

\begin{figure}\centering
\newlength{\minipwidth}\setlength{\minipwidth}{.33\textwidth}
\newlength{\figwidth}\setlength{\figwidth}{.3\textwidth}
\begin{minipage}{65pt}
$L^1$ projection
\end{minipage}\hfill
\begin{minipage}{\minipwidth}\centering
\includegraphics[width=\figwidth]{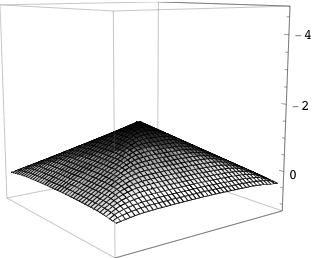}
\end{minipage}\hfill
\begin{minipage}{\minipwidth}\centering
\includegraphics[width=\figwidth]{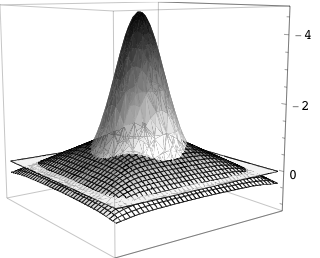}
\end{minipage}\\[10pt]

\begin{minipage}{65pt}
$L^2$ projection
\end{minipage}\hfill
\begin{minipage}{\minipwidth}\centering
\includegraphics[width=\figwidth]{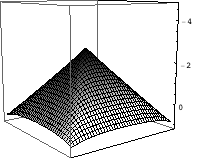}
\end{minipage}\hfill
\begin{minipage}{\minipwidth}\centering
\includegraphics[width=\figwidth]{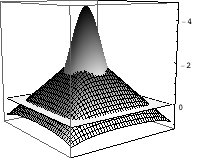}
\end{minipage}\\[10pt]

\begin{minipage}{65pt}
$L^\infty$ projection
\end{minipage}\hfill
\begin{minipage}{\minipwidth}\centering
\includegraphics[width=\figwidth]{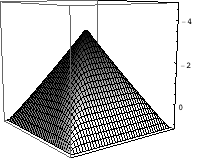}
\end{minipage}\hfill
\begin{minipage}{\minipwidth}\centering
\includegraphics[width=\figwidth]{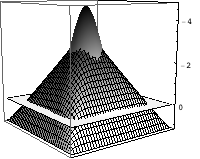}
\end{minipage}\\[10pt]

\begin{minipage}{65pt}
$H^1$ projection
\end{minipage}\hfill
\begin{minipage}{\minipwidth}\centering
\includegraphics[width=\figwidth]{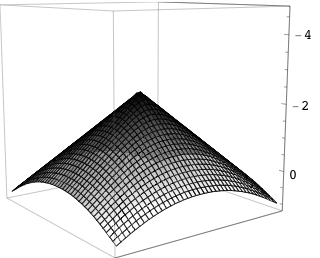}
\end{minipage}\hfill
\begin{minipage}{\minipwidth}\centering
\includegraphics[width=\figwidth]{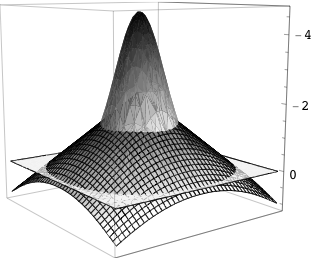}
\end{minipage}\\[10pt]

\begin{minipage}{65pt}
$H_0^1$ projection
\end{minipage}\hfill
\begin{minipage}{\minipwidth}\centering
\includegraphics[width=\figwidth]{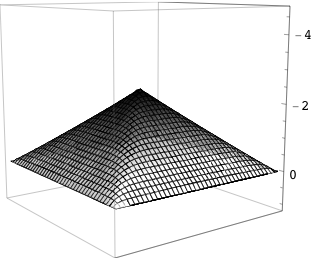}
\end{minipage}\hfill
\begin{minipage}{\minipwidth}\centering
\includegraphics[width=\figwidth]{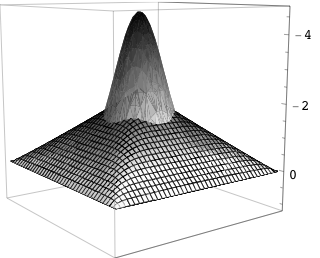}
\end{minipage}

\caption{Projections and its comparisons to the function of figure~\ref{fig:fig2-f}, shown upside down.}
\label{fig:fig2}
\end{figure}

We show the graph of the original function in figure~\ref{fig:fig2-f}, and that of the resulting $L^1$, $L^2$, $L^\infty$, $H^1$ and $H_0^1$ projections in figure~\ref{fig:fig2}.
As in the original article, these graphs are shown upside down.
The interested reader may observe that our results are qualitatively different from those in~\cite[p.~304]{carlier}.

\begin{figure}\centering
\begin{minipage}[b]{.3\textwidth}\centering
\includegraphics[width=.9\textwidth]{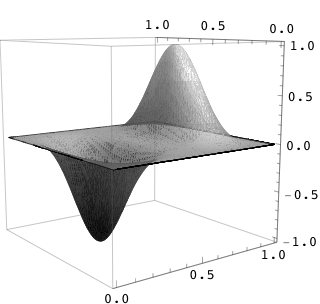}\\ (a)
\end{minipage}
\hfill
\begin{minipage}[b]{.3\textwidth}\centering
\includegraphics[width=.9\textwidth]{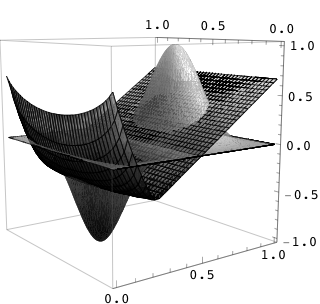}\\ (b)
\end{minipage}
\hfill
\begin{minipage}[b]{.3\textwidth}\centering
\includegraphics[width=.9\textwidth]{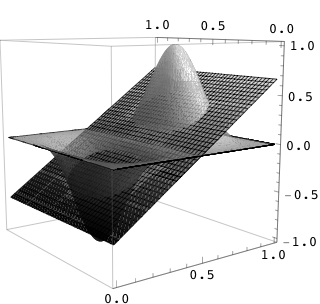}\\ (c)
\end{minipage}
\caption{Non-uniqueness of the $L^\infty$ projection.
(a) original function, (b) and (c) two optimal solutions.}
\label{fig:nonu}
\end{figure}

Using a similar function, but with more symmetries,
\[
   g(x_1,x_2) = -\sin^3 (2\pi x_1)\,\sin^3(\pi x_2)
   \qquad\text{for $(x_1,x_2)\in\Q$},
\]
we show in figure~\ref{fig:nonu} that the $L^\infty$ projection need not be unique.
The solution in~(b) is the one obtained by the semidefinite programming code.

In all these examples, we used a mesh with $41\times 41$ nodes.
The times ranged from about $40$ seconds for the $L^\infty$ projections to about $200$ seconds for the $L^2$ projection.
The $H^1$ and $H_0^1$ projections took about the same time, near $140$ seconds.

\subsection{Fitting data}\label{sec:fit}

As mentioned in the introduction, many problems in science are modelled via convex functions, raising the question of how measured data (usually non convex) can be approximated by a convex function.

Often the fitness to data is done parametrically.
For example, by assuming that the underlying function is a linear combination of some given polynomials, and minimizing over all possible parameters in a convenient norm.
Even in this case, approximating the data by a linear combination which is also convex---but otherwise arbitrary---may be challenging.

In this section we show how our numerical scheme could be used for fitting data on a regular mesh by discrete functions having a positive semidefinite discrete Hessian.
Though the resulting discrete functions may not be extended to a convex function of continuous variables, as shown in example~\ref{example:hessnoconv}, the underlying ``true'' convex function might be well captured by the discrete function.

\begin{figure}\centering
\setlength{\minipwidth}{.45\textwidth}
\setlength{\figwidth}{.4\textwidth}

\begin{minipage}{\minipwidth}\centering
\includegraphics[width=\figwidth]{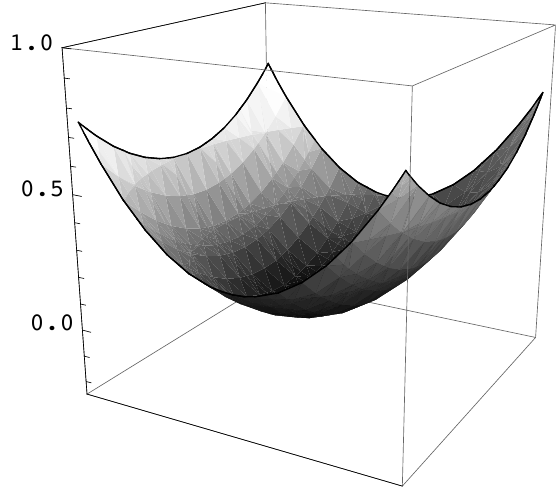}\\ (a)
\end{minipage}\hfill
\begin{minipage}{\minipwidth}\centering
\includegraphics[width=\figwidth]{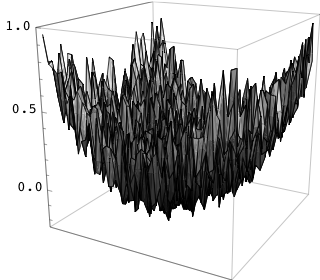}\\ (b)
\end{minipage}\\[10pt]

\begin{minipage}{\minipwidth}\centering
\includegraphics[width=\figwidth]{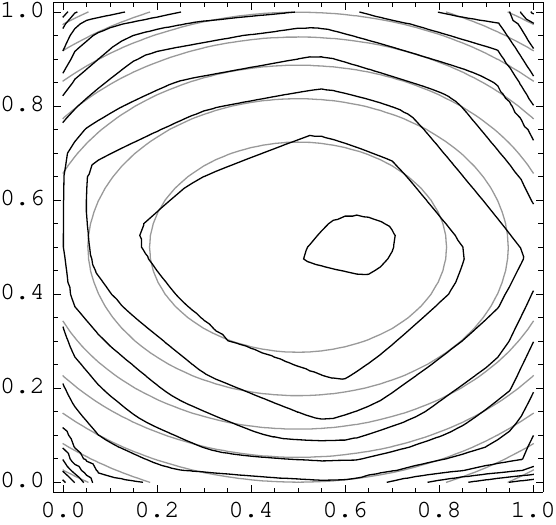}\\ (c)
\end{minipage}\hfill
\begin{minipage}{\minipwidth}\centering
\includegraphics[width=\figwidth]{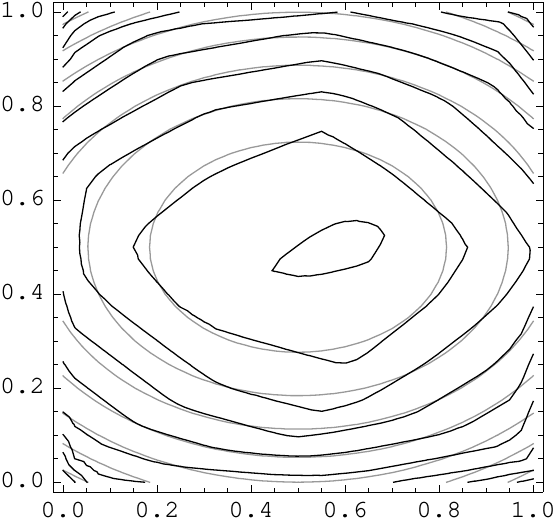}\\ (d)
\end{minipage}\\[10pt]

\begin{minipage}{\minipwidth}\centering
\includegraphics[width=\figwidth]{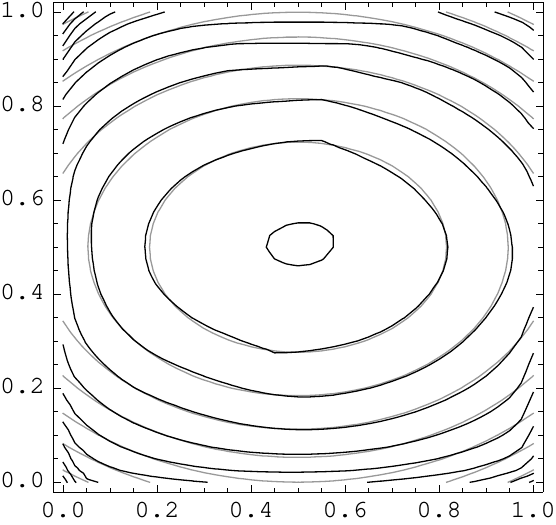}\\ (e)
\end{minipage}\hfill
\begin{minipage}{\minipwidth}\centering
\includegraphics[width=\figwidth]{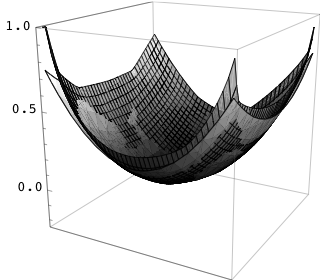}\\ (f)
\end{minipage}

\caption{Fitting a convex surface to perturbed data.
Original function (a),
perturbed data (b),
level curves of the $L^1$ projection (c),
level curves of the $L^2$ projection (d),
level curves of the $L^\infty$ projection (e),
3D graphs of the unperturbed function and $L^\infty$ projection (f).}
\label{fig:noise}
\end{figure}

In the tests we show, we perturbed with random noise the values of the function
\[ f(x_1,x_2) = (x_1 - 1/2)^2 + 2\,(x_2 - 1/2)^2, \]
on a regular mesh of $41\times 41$ nodes.
The original function and the perturbed data are represented in figure~\ref{fig:noise}~(a) and (b), respectively.\footnote{%
   It must be kept in mind that, though the data is discrete, the
   graphic software makes its own interpolations for drawing surfaces
   and level curves.}

We chose to use a uniformly distributed noise between $-\varepsilon$ and $\varepsilon$ for each node, where $\varepsilon = 10\,h$ ($h$ being the mesh size).
In this way we simulate some intrinsic measurement noise whose distribution is presumably known, and that the mesh has been chosen finer than the measurement noise.
Needless to say, our choice of noise is quite arbitrary, and there are many other possibilities to choose from.

In figure~\ref{fig:noise}~(c), (d) and~(e) we show the resulting level curves for the $L^1$, $L^2$ and $L^\infty$ discrete projections (respectively), with the contours of the original unperturbed function in a lighter gray.
The $L^2$ discrete projection may be considered as a variant of a least squares approximation, but the others are not seen as often.

The times spent by the semidefinite code were similar to those of section~\ref{sec:projections:2}, about $170$, $200$ and $40$ seconds (respectively).
As is to be suspected, for the same mesh but smaller noise (e.g., taking $\varepsilon = h$), the results are better and the times smaller.

Comparing the level curves with those of the unperturbed function (figure~\ref{fig:noise}~(c), (d) and~(e)), it is apparent that the $L^\infty$ discrete projection seems to give the closest fit.
Notice that this projection is also faster than the $L^1$ or $L^2$ discrete projections by a factor of $4$ to $5$.

For the run shown, the maximum absolute value of the differences between the perturbed data and the original function on the mesh is about $0.25$ ($10$ times the mesh size, as designed), whereas the maximum absolute value of the differences between the $L^\infty$ projection and the original function  is $0.42$, attained at $(0,1)$.

However, as can be seen in figure~\ref{fig:noise}~(e), the $L^\infty$ projection gives a good fit to the unperturbed function in the interior nodes, and is actually between $\pm 0.018$ (somewhat smaller than the mesh size) for nodes of the square $0.1 \le x_1, x_2\le 0.9$, but deteriorates near the boundary (as do the other projections).
This is perhaps better appreciated in figure~\ref{fig:noise}~(f), where the graph of the $L^\infty$ discrete projection is compared to that of the unperturbed original function.

The ``overshooting'' at the boundary is a typical and known phenomenon when fitting data by convex functions.
See, for instance, the article by Meyer~\cite{Me06}, where the one dimensional case is discussed.

\section*{Acknowledgements}

\begin{itemize}
\item
Our gratitude to A.~Manelli for bringing the monopolist problem~\ref{prob:monopolist} to our attention, and for his many enlightening comments on the subject.

\item
Our thanks to the anonymous referees for their comments which made the presentation more clear, and for suggesting the applications in section~\ref{sec:fit}.
\end{itemize}



{\setlength{\parindent}{0pt}
\small

\medskip
\rule{.5\textwidth}{1.5pt}

\medskip
\textbf{Affiliations}
\begin{description}
\item[N\'{e}stor E. Aguilera:] Consejo Nacional de Investigaciones Cient\'{\i}ficas y T\'{e}cnicas and Universidad Nacional del Litoral, Argentina.\\
e-mail: \url{aguilera@santafe-conicet.gov.ar}

\item[Pedro Morin (corresponding author):] Consejo Nacional de Investigaciones Cient\'{\i}ficas y T\'{e}cnicas and Universidad Nacional del Litoral, Argentina.\\
e-mail: \url{pmorin@santafe-conicet.gov.ar}
\end{description}
}

\end{document}